\documentclass[times,preprint,3p,11pt]{elsarticle}
\usepackage{etoolbox}
\usepackage[title]{appendix}
\usepackage{amscd,amsmath,amssymb,amsthm,mathrsfs,bbm,listings}
\usepackage[hidelinks]{hyperref}
\usepackage{color}
\usepackage{epsfig}
\usepackage{multirow}
\usepackage{verbatim}
\usepackage{subfig}
\usepackage{soul}
\usepackage{pgfplots}
\usepackage{cuted}
\pgfplotsset{compat=1.11}

\numberwithin{equation}{section}

\allowdisplaybreaks
\usepackage[displaymath,mathlines,pagewise]{lineno}
\usepackage[ruled,vlined,boxed,linesnumbered,longend]{algorithm2e}

\hypersetup{
    colorlinks=true,
    pdfborder={0 0 0},
    allcolors =blue,
}

\newcommand{\Diff}[2]{\mathcal{D}_{#1}\left(#2\right)}
\newcommand{\React}[2]{\mathcal{R}_{#1}\left(#2\right)}
\newcommand{\dof}[2]{\text{dof}_{#1}\left(#2\right)}
\newcommand{\diam}[1]{\mbox{diam}\left(#1\right)}

\newcommand{\Vke}[2]{\widetilde{V}_{#1}\left(#2\right)} 
\newcommand{\VkS}[2]{V_{#1}^S\left(#2\right)}
\newcommand{\piB}[2]{\pi^{\partial}_{#1} \left(#2\right)}
\newcommand{\ppiB}[1]{\pi^{\partial}_{#1}}
\newcommand{\piO}[2]{\pi^{0}_{#1} \left(#2\right)}
\newcommand{\ppiO}[1]{\pi^{0}_{#1} }
\newcommand{\piN}[2]{\pi^{\nabla}_{#1} \left(#2\right)}
\newcommand{\piS}[2]{\pi^{S}_{#1} \left(#2\right)}
\newcommand{\ppiN}[1]{\pi^{\nabla}_{#1}}
\newcommand{\Ph}[1]{\mathcal{P}_{h}\left(#1\right)}
\newcommand{\Ih}[1]{\mathcal{I}^{#1}_{h}}
\newcommand{\ac}[2]{a\left(#1;#2\right)}
\def\UhE{\uu{U}_{h, E}}  
\newcommand{\acE}[2]{a^E\left(#1;#2\right)}
\newcommand{\ah}[2]{a_h\left(#1;#2\right)}
\newcommand{\ahE}[2]{a_h^E\left(#1;#2\right)}
\newcommand{\saE}[2]{s_a^E\left(#1;#2\right)}
\newcommand{\mc}[2]{m\left(#1;#2\right)}
\newcommand{\mcE}[2]{m^E\left(#1;#2\right)}
\newcommand{\mh}[2]{m_h\left(#1;#2\right)}
\newcommand{\mhE}[2]{m_h^E\left(#1;#2\right)}
\newcommand{\smE}[2]{s_m^E\left(#1;#2\right)}
\newcommand{\Seminorm}[2]{\left| #1 \right|_{#2}}
\newcommand{\Norm}[2]{\left\| #1 \right\|_{#2}}
\newcommand{\EFC}[2]{\mathcal{C}^{#1}\left(#2\right)}
\newcommand{\dfdt}[1]{\frac{d#1}{dt}} 
\newcommand{\dfhdt}[1]{\frac{\partial #1}{\partial t}}
\newcommand{\Pp}[2]{\IP_{#1}\left(#2\right)}
\newcommand{\DotProd}[2]{{\left(#1,#2\right)}}

\newcommand{\uu}[1]{\hbox{\boldmath$#1$}}

\def\bM{\mathbf{M}}

\def\bI{\mathbf{I}}
\def\bA{\mathbf{A}}

\def\PPh{\mathcal{P}_{h}}
\def\aah{a_h}
\def\mmh{m_h}
\def\uh{u_h} 
\def\Uh{\uu{U}_h} 
\def\dofi{\mbox{dof}_i}
\def\dofj{\mbox{dof}_j}
\def\bxi{\uu{\xi}}
\def\bx{\uu{x}}
\def\ETAL{{\em et. al.,\ }}
\def\IR{\mathbb R}

\def\IP{\mathbb P}
\def\Th{{\mathcal T}_h}

\def\vh{v_h}
\def\dx{d\uu{x}}
\def\dS{dS}
\def\QED{\hfill\ensuremath{\square}}
\def\nv{\vec{\uu{n}}}
\def\PROOF{{\em Proof. }}

\newtheorem{assumption}{Assumption}
\newtheorem{proposition}{Proposition}

\newtheorem{theorem}{Theorem}
\newtheorem{remark}{Remark}
\newtheorem{lemma}{Lemma}
\newcommand{\ORDER}[1]{{\mathcal O}\left(#1\right)}

\begin{document}
\begin{frontmatter}
\title {High-order interpolatory Serendipity Virtual Element Method for 
semilinear parabolic equations.}

\author[UNIPV,USI]{Sergio A.~G\'omez  \corref{corrauthor}}
\ead{sergio.gomez01@universitadipavia.it,gomezs@usi.ch}

\address[UNIPV]{
Dipartimento di Matematica ``F. Casorati'', Universit\`a di Pavia,\\
 Via Ferrata 5, 27100, Pavia, Italy.}
 \address[USI]{
 Faculty of Informatics, Universit\`a della Svizzera italiana, \\
 Via Buffi 13, 6900, Lugano, Switzerland.} 
\cortext[corrauthor]{Corresponding author}
\begin{linenomath}
\begin{abstract}
We propose an efficient method for the numerical approximation of a general 
class 
of two dimensional semilinear parabolic problems on polygonal meshes. The 
proposed approach takes advantage of the properties of the serendipity version 
of the Virtual 
Element Method, which not only reduces the number of 
degrees of freedom compared to the original Virtual Element Method, {but 
also allows the introduction of an approximation of the 
nonlinear term that is computable from the degrees of freedom of the 
discrete solution with a low computational cost}, thus significantly 
improving the efficiency of the method. An error analysis for 
the semi-discrete formulation is carried out, and an optimal 
estimate for the error in the $L_2$-norm is obtained. The accuracy and 
efficiency of the proposed method when combined with a second order Strang 
operator splitting time discretization is illustrated in our numerical 
experiments, with approximations up to order $6$.
\end{abstract}
\end{linenomath}

\end{frontmatter}
\section{Introduction \label{SECT::INTRODUCTION}}
In this work we present an interpolatory or quasi-interpolatory Serendipity Virtual Element Method 
(S-VEM) applied to semilinear parabolic equations on a space--time 
domain $Q_T = \Omega \times (0, 
T)$, 
where $\Omega \subset \IR^2$ is a polygonal domain and $T > 0$
\begin{linenomath}
\begin{subequations}
\label{EQN::MODEL-PROBLEM}
\begin{align}
\dfhdt{u} - \Delta u + f(u) & = 0, \ \ \qquad 
\mbox{ in } Q_T,\\
\nabla u \cdot \nv &= 0, \quad\quad\ \ \mbox{ on } \partial \Omega \times (0, 
T],\\
u(\bx, 0) & = u_0(\bx), \quad \mbox{ in }\Omega.
\end{align}
\end{subequations}
\end{linenomath}
{The nonlinear function $f: \IR \rightarrow \IR$ is assumed to be globally 
Lipschitz continuous, i.e., there exists a constant $L_f > 0$ such that the 
following 
bound holds
\begin{equation}
\label{EQN::LIPCHITZ-CONDITION}
\left|f(x) - f(y)\right|\ \leq\ L_f\left|x- y \right| \quad   
\forall x, y \in 
\IR.
\end{equation}
}
The model \eqref{EQN::MODEL-PROBLEM} is found in many important applications 
such as: battery modeling \cite{Tenno_Tenno_Suntio_2001}, crystals growth 
\cite{Kobayashi_1993}, population dynamics \cite{Neubert_Caswell_2000}, and in 
many other models in chemistry 
\cite{Turing_1990,Mikhailov_Showalter_2006} and biology 
\cite{Marcon_Sharpe_2012}. However, given the different nature of nonlinear 
terms, the task of finding exact solutions for such 
kind of problems becomes extremely demanding or even impossible. For that 
reason, there is a high interest in the development of efficient, accurate and 
robust numerical methods to approximate their solution. Since this work 
specifically concerns to the advantages of the serendipity version of the 
Virtual Element Method applied to the problem \eqref{EQN::MODEL-PROBLEM}, an 
extensive list 
of numerical methods previously applied to this problem is out of our scope.

The Virtual Element Method (VEM) is a novel technique for the numerical 
approximation of PDEs, introduced by Beir\~ao da Veiga \ETAL in 
\cite{Beirao_Brezzi_Cangiani_Manzini_Marini_Russo_2013}, for an elliptic 
problem, and 
can be seen as a sensible extension of the classical finite element method to 
meshes with almost general polygonal elements. Discrete VEM spaces 
contain non-polynomial functions; however, such functions are not needed to be 
expressly known, as the discrete operators are computed through 
projections onto the space of piecewise polynomials of a given degree, which 
are 
computable using only some suitably chosen degrees of freedom (DoFs). Besides 
the 
advantages 
that come from the versatility of polygonal meshes, such as the 
natural use of ``non-conforming'' grids, more efficient and easier adaptivity 
and geometric approximation, robustness to mesh 
deformation, among others; the Virtual Element Method  also allows for the 
imposition of conformity conditions on the global discrete spaces without 
struggling to explicitly compute their basis functions.

So far, the Virtual Element Method has been successfully applied to 
many important physical applications. In particular, recent efforts have been 
devoted to show the accuracy and advantages of this method in the numerical 
approximation of the solution to nonlinear problems such as: the Cahn--Hilliard 
equation 
\cite{Antonietti_DaVeiga_Scacchi_Verani_2016,Liu_He_Chen_2020}, 
models in cardiology 
\cite{Anaya_Bendahmane_Mora_Sepulveda_2020}, nonlinear elasticity  
\cite{Bellis_Wriggers_Hudobivnik_2019}, the nonlinear Brinkman equation 
\cite{Gatica_Munar_Sequeira_2018,Munar_Sequeira_2020}, bulk-surface 
reaction--diffusion 
systems \cite{Frittelli_Madzvamuse_Sgura_2021}, pattern formation 
\cite{Dehghan_Gharibi_2021}; and semilinear 
elliptic 
\cite{Cangiani_Chatzipantelidis_Diwan_Geourgoulis_2020,
Adak_Natarajan_Natarajan_2019},
 hyperbolic \cite{Adak_Natarajan_Kumar_2019_b} and parabolic 
\cite{Adak_Natarajan_Kumar_2019} equations. 

In this work we aim to extend the idea of Adak and Natarajan in 
\cite{Adak_Natarajan_2020} to high-order approximations.  In 
\cite{Adak_Natarajan_2020}, the authors proposed a VEM discretization 
for the sine--Gordon equation with an interpolatory approximation of 
the nonlinear term, thus significantly reducing the computational cost of the 
method at 
each time step. Nevertheless, the main limitation of the technique in 
\cite{Adak_Natarajan_2020} is that it is only valid for 
approximations with $k = 1$, i.e., with the same order of convergence as 
polynomial approximations of degree $k = 1$. This is due to the fact that for 
$k \geq 2$, the method requires some internal-moment  DoFs, which 
unfortunately prevents {a direct extension to high order approximations} 
(see 
Remark 3.3 
and Section 7 in \cite{Adak_Natarajan_2020} dedicated to discuss this 
limitation). Our idea to overcome this severe restriction relies on the 
serendipity variation of the VEM, introduced by Beir\~ao da Veiga \ETAL\ in 
\cite{DaVeiga_Brezzi_Marini_Russo_2016}, and later discussed by Russo in 
\cite{Russo_2016}. The main motivation of the S-VEM is indeed to reduce the 
number 
of internal-moment DoFs. {Moreover, under certain conditions on the mesh, it 
is possible} to completely eliminate them. It is also worth mentioning that 
the Serendipity VEM on quadrilaterals does not suffer from distortion as it is 
common for the serendipity FEM, see 
\cite{DaVeiga_Brezzi_Marini_Russo_2016}.

The main novelty and features of the proposed scheme are summarized as 
follows:
\begin{enumerate}
	\item[a)] To the best of our knowledge, this is the first time to use 
the S-VEM 
as 
spatial discretization for semilinear parabolic problems; for which the 
enhanced VEM has been preferred.
	\item[b)] An interpolant of the nonlinear term in the S-VEM space is 
introduced in the semi-discrete formulation. Under a {condition} on the 
degree {of accuracy $k$ that is associated to the geometric properties of the 
mesh}, such interpolant is computed {by simply evaluating the nonlinear 
term $f(\cdot)$ at the DoFs of the discrete solution.  Thus significantly 
reducing  the  computational cost of the method, as it completely avoids} the 
use of quadratures at each time step.
	\item[c)] An optimal error estimate for the semi-discrete formulation in 
the 
$L_2$-norm is proven in spite of the use of such {interpolant} to 
approximate the nonlinear term.
	\item[d)] When the time variable is discretized by the symmetric Strang - 
operator splitting (SS-OS) time marching scheme, the nonlinear substeps can be 
decomposed as a set of completely independent one dimensional nonlinear 
problems, {which makes} the method naturally suitable for parallel 
implementations.
	\item[e)] {In those elements of the mesh (if any) where the condition on 
the degree $k$ is not satisfied, the interpolant of the nonlinear term is not 
computable from the DoFs of the discrete solution. In that case, we use a quasi-interpolatory approximation of the nonlinear term that also belongs to the local 
VEM space but is computable. Optimal error 
estimates and suitability for a parallel 
implementation are preserved.}
\end{enumerate}

The paper is structured as follows: in Section \ref{SECT::SERENDIPITY-VEM} we 
present the basic ideas and necessary projections for the description of the  
proposed 
method {in the case when no internal-moment DoFs are needed, that throughout 
the paper we will refer to as ``the ideal case''}. Optimal error estimates of 
order $\mathcal{O}(h^{k+1})$ in the 
$L_2$-norm are proven 
for the S-VEM semi-discrete formulation in Section \ref{SECT::ERROR-ANALYSIS}. 
In 
Section \ref{SECT::FULLY-DISCRETE}, an efficient fully-discrete scheme, 
obtained by combining our interpolatory S-VEM discretization in space with a 
symmetric Strang - Operator Splitting time marching scheme is presented. 
{The 
extension of the method to the general case when some internal-moment DoFs are 
needed, as well as the most important differences in the error estimate and the 
fully-discrete scheme are discussed in Section 
\ref{SECT::HIGH-ORDER}.}
Some numerical experiments, validating the accuracy and efficiency of the 
proposed method are included in Section \ref{SECT::NUMERICAL-EXPERIMENTS}. 
We end this work with some concluding remarks in Section 
\ref{SECT::CONCLUSIONS}.

\section{Serendipity VEM discretization\label{SECT::SERENDIPITY-VEM}}
Let $\Th$ be a polygonal partition of $\Omega$ and let 
$h := \max \left\{h_E \ | \ E \in \Th\right\}$ be the mesh size, where $h_E$ 
denotes the diameter of $E$. 
We first define, for each polygon 
$E\in \Th$, the following 
enlarged local Virtual Element space 
\cite{Ahmad_Alsaedi_Brezzi_Marini_Russo_2013}:
\begin{equation*}
\label{EQN::INFLATED-VEM}
	\Vke{k}{E} := \left\{v \in \EFC{0}{\overline{E}}: \ v|_e 
\in \Pp{k}{e}\ \forall \mbox{ edge } e, \ \Delta v \in \Pp{k}{E}\right\},
\end{equation*}
where $\Pp{k}{e}$ and $\Pp{k}{E}$ denote the spaces of polynomials of degree at 
most $k \geq 1$ on $e$ and $E$, respectively.

The DoFs uniquely identifying a function $v \in \Vke{k}{E}$ 
are choosen as {the following linear functionals}
\begin{enumerate}
	\item[i)] The values of $v$ at the vertices of $E$.
	\item[ii)] The values of $v$ at {the $(k - 1)$ internal Gauss-Lobatto nodes 
on each edge $e$ of $E$}.
	\item[iii)] The internal-moments:\  $\displaystyle \frac{1}{|E|}\int_E v 
{m_\alpha^E \dx, \ \alpha = 1, \ldots, r_k}$,

{where $\left\{m_\alpha \right\}_{\alpha = 1}^{r_k}$ is a basis of 
$\Pp{k}{E}$ and $r_k := \dim(\Pp{k}{E})$.}
\end{enumerate}

The space 
$\Vke{k}{E}$ requires many more internal DoFs than the 
original 
VEM space presented in \cite{Beirao_Brezzi_Cangiani_Manzini_Marini_Russo_2013}, 
but it readily provides enough information to 
compute the 
$L_2$-projection of any $v\in \Vke{k}{E}$ onto $\Pp{k}{E}$. In practice, a 
subspace of $\Vke{k}{E}$ still containing all polynomials of degree at most 
$k$ on $E$, is used as local VEM 
space. 
The basic idea to construct such subspace is to take the set of functions in 
$\Vke{k}{E}$ sharing some internal-moment DoFs with their projection
onto the space $\Pp{k}{E}$; which gives origin to the so called, 
\textit{enhanced} \cite{Ahmad_Alsaedi_Brezzi_Marini_Russo_2013} and 
\textit{serendipity} \cite{DaVeiga_Brezzi_Marini_Russo_2016} versions of {the} 
VEM.

{We first focus on the ideal case, where the Serendipity VEM does not 
require 
any internal-moment DoFs.}
An integer number $\eta_E \geq 3$ is associated to each element $E\in \Th$, 
where $\eta_E$ is defined as 
the 
number of distinct straight lines containing at least one edge of $E$. In 
particular, if $E$ is an $N$-sided strictly convex polygon without split edges, 
then
$\eta_E= N$. {In the ideal case, the degree of accuracy $k$ satisfies the 
condition: $k < \min \left\{\eta_E \ 
|\ E \in \Th\right\}$}; which in the spirit of {the} Serendipity VEM, allows 
{the definition of}
a global discrete space whose associated DoFs are all node evaluations at the 
skeleton of the mesh, without requiring any internal-moment degree of freedom 
from the set iii). 

The following projectors are needed to define the local S-VEM space 
and to present our semi-discrete formulation:
\begin{itemize}
	\item[$\bullet$] The Ritz--Galerkin projection $\ppiN{k, E} : {H^1(E)}
\rightarrow \Pp{k}{E}$ defined as follows
\begin{align*}
\label{EQN::PI_NABLA_PROJECTION-1}
\int_E \nabla \left(\piN{k, E}{v} - v\right) \cdot \nabla p_k \dx = 0 \quad 
\forall p_k \in \Pp{k}{E}, \\
\oint_{\partial E} \piN{k, E}{v} \dS = \oint_{\partial E} v \dS, \ (k = 
1), \ \mbox{ or } \
\int_E \piN{k, E}{v} \dx & = \int_E v \dx, \ (k > 1).
\end{align*}
Using the Green's formula, 
{the} projection $\piN{k, E}{\cdot}$ {is computable for any $v \in 
\Vke{k}{E}$} using the degrees 
of freedom i), ii) and iii), see 
\cite[Sect. 4.5]{Beirao_Brezzi_Cangiani_Manzini_Marini_Russo_2013}.
	\item[$\bullet$] The standard $L_2$-orthogonal projector $\ppiO{k, E}: 
{L_2(E)}
\rightarrow \Pp{k}{E}$, defined by
\begin{equation*}
\label{EQN::L2-PROJECTION}
\int_E \left(\piO{k, E}{v} - v\right)p_k \dx = 0 \quad \forall p_k \in 
\Pp{k}{E},
\end{equation*}
which is directly computable from the set of DoFs iii).

	\item[$\bullet$] The ``boundary'' projector $\ppiB{k, E} : {H^1(E)}
\rightarrow 
\Pp{k}{E}$ such 
that
\begin{equation}
	\label{EQN::BNDRY-PROJECTION}
	\oint_{\partial E} \left(\piB{k, E}{v} - v\right)p_k\ \dS = 0 
\quad \forall p_k \in \Pp{k}{E},
\end{equation}
that is well-defined under the condition $k < \eta_E$, and can be computed 
using only the boundary DoFs in sets i) and ii).

\end{itemize}
{In the ideal case,} the local S-VEM space $\VkS{k}{E} \subset \Vke{k}{E}$ 
is defined as
\begin{equation*}
\label{EQN::SERENDIPITY-VEM-SPACE}
	\VkS{k}{E} := \left\{v \in \Vke{k}{E} : \int_E \left(v - 
\piB{k, E}{v}\right) {m_\alpha^E} \dx \ =\ 0, \ {\alpha = 1, \ldots, 
r_k}\right\},
\end{equation*}
{that only requires the boundary DoFs from sets i) and ii)}. Therefore, for 
an $N$-sided polygon $E$, $\dim\left({\VkS{k}{E}}\right) = 
Nk$. The global S-VEM space is consequently defined as 
\[\VkS{k}{\Th} := 
\left\{\vh \in \EFC{0}{\overline{\Omega}} : \ \vh|_E \in 
\VkS{k}{E} \ \forall E \in \Th\right\}.\]

\begin{figure}[!ht]
\centering
\begin{minipage}{0.32\textwidth}
\begin{center}
\begin{tikzpicture}[scale = .9]
\draw[black, very thick, -] (0,0) -- (2, 4) -- (4,0) -- (0,0);
\foreach \point in {(0,0),(2,4),(4,0),(1,2),(2, 0),(3, 2)}{
    \draw[fill = blue] \point circle (5pt);
}
\node [draw, thick, shape=rectangle, minimum width=1mm, minimum height=1mm, 
fill=red, anchor=center] at (2,1.5) {};
\end{tikzpicture}
\end{center}
\end{minipage}
\begin{minipage}{0.32\textwidth}
\begin{center}
\begin{tikzpicture}[scale=.9]
\draw[black, very thick, -] (0.00, 0.00) -- (0.17, 1.33) -- (0.33, 2.67) -- 
(0.50, 4.00) -- (0.50, 4.00) -- (1.33, 4.00) -- (2.17, 4.00) -- (3.00, 4.00) -- 
(3.00, 4.00) -- (3.67, 2.67) -- (4.33, 1.33) -- (5.00, 0.00) -- (5.00, 0.00) -- 
(3.33, 0.00) -- (1.67, 0.00) -- (0.00, 0.00);
\foreach \point in {(0.00, 0.00), (0.17, 1.33), (0.33, 2.67), (0.50, 4.00), 
(0.50, 4.00), (1.33, 4.00), (2.17, 4.00), (3.00, 4.00), (3.00, 4.00), (3.67, 
2.67), (4.33, 1.33), (5.00, 0.00), (5.00, 0.00), (3.33, 0.00), (1.67, 0.00), 
(0.00, 0.00)
}{
    \draw[fill = blue] \point circle (5pt);
}
\node [draw, thick, shape=rectangle, minimum width=1mm, minimum height=1mm, 
fill=red, anchor=center] at (1.67, 1.37) {};
\node [draw, thick, shape=rectangle, minimum width=1mm, minimum height=1mm, 
fill=red, anchor=center] at (2.07, 2.29) {};
\node [draw, thick, shape=rectangle, minimum width=1mm, minimum height=1mm, 
fill=red, anchor=center] at (2.82, 1.37) {};
\end{tikzpicture}
\end{center}
\end{minipage}
\begin{minipage}{0.32\textwidth}
\begin{center}
\begin{tikzpicture}[scale=.9]
\draw[black, very thick, -] (0.60, 0.07)-- (0.38, 0.74)-- (0.17, 1.40)-- 
(-0.05, 2.07)-- (-0.27, 2.74)-- (-0.27, 2.74)-- (0.30, 3.15)-- (0.87, 3.56)-- 
(1.43, 3.97)-- (2.00, 4.38)-- (2.00, 4.38)-- (2.57, 3.97)-- (3.13, 3.56)-- 
(3.70, 3.15)-- (4.27, 2.74)-- (4.27, 2.74)-- (4.05, 2.07)-- (3.83, 1.40)-- 
(3.62, 0.74)-- (3.40, 0.07)-- (3.40, 0.07)-- (2.70, 0.07)-- (2.00, 0.07)-- 
(1.30, 0.07)-- (0.60, 0.07);
\foreach \point in {(0.60, 0.07), (0.38, 0.74), (0.17, 1.40), (-0.05, 2.07), 
(-0.27, 2.74), (-0.27, 2.74), (0.30, 3.15), (0.87, 3.56), (1.43, 3.97), (2.00, 
4.38), (2.00, 4.38), (2.57, 3.97), (3.13, 3.56), (3.70, 3.15), (4.27, 2.74), 
(4.27, 2.74), (4.05, 2.07), (3.83, 1.40), (3.62, 0.74), (3.40, 0.07), (3.40, 
0.07), (2.70, 0.07), (2.00, 0.07), (1.30, 0.07), (0.60, 0.07)
}{
    \draw[fill = blue] \point circle (5pt);
}
\foreach \point in {(1.25, 1.57),(1.62, 2.22),(2.00, 2.87),(2.00, 2.87),(2.38, 
2.22),(2.75, 1.57),(2.75, 1.57),(2.00, 1.57),(1.25, 1.57)}
{
\node [draw, thick, shape=rectangle, minimum width=1mm, minimum height=1mm, 
fill=red, anchor=center] at \point {};
}
\end{tikzpicture}
\end{center}
\end{minipage}
\caption{Degrees of freedom for  {the original VEM} space with degree $k = 
2$ 
(triangles), 
$k = 3$ (quadrilaterals), and $k = 4$ (pentagons). Blue dots 
represent nodal evaluations, and red squares represent internal moments. 
\label{FIG::DOF}}
\end{figure}
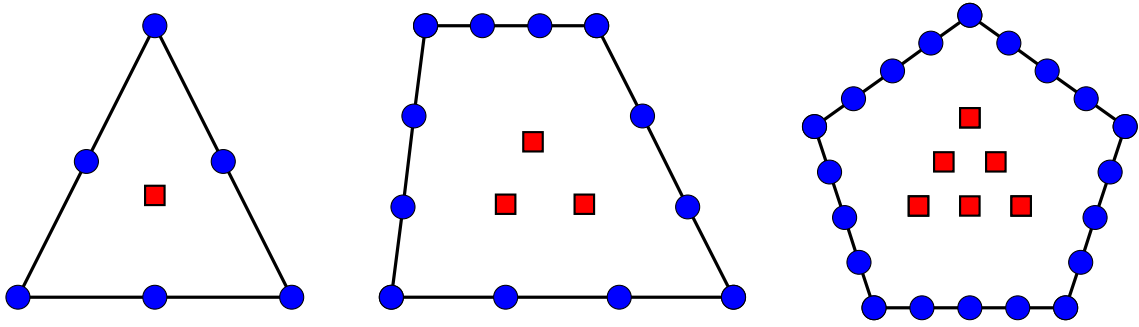

A representation of the DoFs for the {local space of the original VEM in 
\cite{Beirao_Brezzi_Cangiani_Manzini_Marini_Russo_2013}}, for 
different 
$N$-sided polygons and the maximum degree $k$ satisfying the aforementioned 
condition is presented in Fig. \ref{FIG::DOF}. {We emphasize that, in all 
these cases, the local space $\VkS{k}{E}$ does not require any of the 
internal-moment DoFs represented by red squares in the figure.
}

As mentioned before, any function $v \in \VkS{k}{E}$ is uniquely determined by 
its boundary DoFs. Denoting by $d_{k,E}^S 
= \dim(\VkS{k}{E})$, and numbering the nodes associated to the DoFs as\ 
$\bxi_i$, with $\ i = 1, \ldots, d_{k,E}^S$, we 
can define the linear functionals $\dofi : \VkS{k}{E} \rightarrow \IR$ as 
\begin{equation}
\label{EQN::DOF-DEF}
\dofi(v) := v(\bxi_i)  \quad \forall v \in \VkS{k}{E}.
\end{equation}
A natural basis arise, by taking the canonical basis functions 
$\left\{\phi_i\right\}_{i = 1}^{d_{k,E}^S}$
satisfying 
\begin{equation*}
\label{EQN::INTERPOLATORY-BASIS}
\dofi(\phi_j) = \delta_{ij}, \quad i, j = 1, \ldots, 
d_{k, E}^S.
\end{equation*}
The following interpolatory representation is then obtained for each $v \in 
\VkS{k}{E}$
\begin{equation}
\label{EQN::DOFI}
v(\bx) = \sum_{i = 1}^{d_{k,E}^S} \dofi(v) \phi_i(\bx) = \sum_{i = 
1}^{d_{k,E}^S} v(\bxi_i) \phi_i(\bx);
\end{equation}
such representation, in turn, allows {us} to define the {interpolant} 
operator 
$\Ih{k} : \EFC{0}{\overline{\Omega}} \rightarrow \VkS{k}{\Th}$ whose 
restriction 
to each element $E \in \Th$ is {defined as follows: for all $g \in 
\EFC{0}{\overline{E}}$,}
\begin{equation}
	\Ih{k, E} g(\bx) := {\sum_{i = 1}^{d_{k, E}^S} \dof{i}{g} \phi_i(\bx)} = 
\sum_{i = 1}^{d_{k, E}^S} g(\bxi_i) \phi_i(\bx).
	\label{EQN::INTERPOLATION-OPERATOR}
\end{equation}
\subsection{Semi-discrete formulation}
The weak formulation of the model problem \eqref{EQN::MODEL-PROBLEM} is: find 
$u\in L_2\left(0, T, H^1(\Omega)\right)$ with $u_t \in L_2\left(0, T, 
H^{-1}(\Omega)\right)$ such that
\begin{equation}
\label{EQN::WEAK-FORMULATION}
\mc{\dfhdt{u}}{v} + \ac{u}{v} + \mc{f(u)}{v} = 0 \quad \forall v \in 
H^1\left(\Omega\right),
\end{equation}
where $\ac{\cdot}{\cdot} : H^1(\Omega) \times H^1(\Omega) \rightarrow \IR$ 
and $\mc{\cdot}{\cdot} : L_2(\Omega) \times L_2(\Omega) \rightarrow \IR$ 
are the bilinear forms defined as
\[\ac{u}{v} := \int_\Omega \nabla u 
\cdot \nabla v \dx, \quad \mc{u}{v} := \int_\Omega u v \dx. \]
Analogously, our semi-discrete interpolatory S-VEM formulation seeks an 
approximation $\uh 
\in \VkS{k}{\Th}$ such that for 
all test functions $\vh \in \VkS{k}{\Th}$ it satisfies 
\begin{linenomath}
\begin{subequations}
\label{EQN::VEM-SEMIDISCRETE-FORMULATION}
\begin{align}
	\label{EQN::VEM-SEMIDISCRETE}
	\mh{\dfhdt{\uh}}{\vh} + \ah{\uh}{\vh} + \mh{\Ih{k} f(\uh)}{\vh} 
& = 0, \\
	\label{EQN::VEM-INITIAL-CONDITION}
	\uh^0 & = \Ih{k} u_0,
\end{align}
\end{subequations}
\end{linenomath}
where the bilinear forms $\aah:\VkS{k}{\Th} \times \VkS{k}{\Th} 
\rightarrow \IR$ and $\mmh: \VkS{k}{\Th} \times \VkS{k}{\Th} 
\rightarrow \IR$ are constructed as the sum of local contributions as
\begin{linenomath}
\begin{equation*}
	\ah{\uh}{\vh} = \sum_{E \in \Th} \ahE{\uh}{\vh}, \quad \quad \mh{\uh}{\vh} 
= \sum_{E \in \Th} \mhE{\uh}{\vh}.
\end{equation*}
\end{linenomath}
For each element $E\in \Th$, the restrictions $\ahE{\cdot}{\cdot}$ and 
$\mhE{\cdot}{\cdot}$ are split 
into a consistency and a stability parts by
\begin{linenomath}
\begin{subequations}
\begin{align}
	\label{EQN::AH-BILINEAR}
	\ahE{\uh}{\vh} := \acE{\piN{k, E}{\uh}}{\piN{k, E}{\vh}} + \saE{\left(I - 
\ppiN{k, E}\right) \uh}{\left(I - 
\ppiN{k, E}\right) \vh}, \\
	\label{EQN::MH-BILINEAR}
	\mhE{\uh}{\vh} := \mcE{\piO{k, E}{\uh}}{\piO{k, E}{\vh}} + \smE{\left(I - 
\ppiO{k, E}\right) \uh}{\left(I - 
\ppiO{k, E}\right) \vh},
\end{align}
\end{subequations}
\end{linenomath}
where $I$ denotes the identity operator and the stabilization terms 
$\saE{\cdot}{\cdot}$, $\smE{\cdot}{\cdot}$ are symmetric bilinear forms 
scaling as $\acE{\cdot}{\cdot}$ and $\mcE{\cdot}{\cdot}$, respectively; 
more precisely, there exist positive constants $\alpha_1, \alpha_2, \beta_1, 
\beta_2$ such that
\begin{align*}
\alpha_1 \acE{\vh}{\vh} \leq \saE{\vh}{\vh} \leq \alpha_2 \acE{\vh}{\vh} & 
\quad \forall \vh 
\in \VkS{k}{E} \cap \text{Ker}\left(\ppiN{k}\right), \\
\beta_1 \mcE{\vh}{\vh} \leq \smE{\vh}{\vh} \leq \beta_2 \mcE{\vh}{\vh} & \quad 
\forall \vh 
\in \VkS{k}{E} \cap \text{Ker}\left(\ppiO{k}\right).
\end{align*}
In 
fact, there are many possible choices for the stability terms; however, {in our 
implementation} we will 
limit ourselves to use a very simple stabilization proposed in 
\cite{Beirao_Brezzi_Cangiani_Manzini_Marini_Russo_2013}, {namely, the 
\textit{dofi-dofi}}. For a thorough study {of} different stability choices, 
see \cite{DaVeiga_Lovadina_Russo_2017,Mascotto_2018}.
By construction, both $\ah{\cdot}{\cdot}$ and $\mh{\cdot}{\cdot}$ satisfy the 
following two important conditions: 
\begin{itemize}
	\item \textbf{$\bf{k}$-Polynomial consistency:} For every element $E \in 
\Th$ 
we have
	\begin{linenomath}
	\begin{subequations}
	\label{EQN::CONSISTENCY-BILINEAR}
	\begin{align}
	\label{EQN::POLYNOMIAL-CONSISTENCY-A}
		\ahE{\vh}{p_k} = \acE{\vh}{p_k} \quad \forall \vh \in \VkS{k}{E}, \ 
\forall p_k \in \Pp{k}{E},\\
	\label{EQN::POLYNOMIAL-CONSISTENCY-M}
		\mhE{\vh}{p_k} = \mcE{\vh}{p_k} \quad \forall \vh \in \VkS{k}{E}, \ 
\forall p_k \in \Pp{k}{E}.
	\end{align}
	\end{subequations}
	\end{linenomath}
	\item \textbf{Stability:} There exist mesh-independent positive constants 
$\alpha_*, 
\ \alpha^*,\ \beta_*,\ \beta^*$ such that
	\begin{linenomath}
	\begin{subequations}
	\label{EQN::STABILITY-BILINEAR}
	\begin{align}
	\label{EQN::STABILITY-A}
		\alpha_* \acE{\vh}{\vh} \leq \ahE{\vh}{\vh} \leq \alpha^* 
\acE{\vh}{\vh} \quad \forall \vh \in \VkS{k}{E},\\
	\label{EQN::STABILITY-M}
		\beta_* \mcE{\vh}{\vh} \leq \mhE{\vh}{\vh} \leq \beta^* 
\mcE{\vh}{\vh} \quad \forall \vh \in \VkS{k}{E}.
	\end{align}
	\end{subequations}
	\end{linenomath}
\end{itemize}
The last term {in the semi-discrete variational formulation} 
\eqref{EQN::VEM-SEMIDISCRETE} 
satisfies the following {crucial identity}
\begin{linenomath}
\begin{align}
\nonumber
	\mh{\Ih{k} f(\uh)}{\phi_i} 
=  \sum_{E \in \Th} \mhE{\Ih{k, E}{f(\uh)}}{\phi_i} & 
\stackrel{\eqref{EQN::DOFI}}{=} \sum_{E \in \Th} \sum_{j = 
1}^{d_{k, E}^S} \dofj\left(f\left(\uh\right)\right)  \mhE{\phi_j}{\phi_i} 
\nonumber \\
& \stackrel{\eqref{EQN::INTERPOLATION-OPERATOR}}{=} \sum_{E \in \Th} \sum_{j = 
1}^{d_{k, E}^S} f\left(\uh\left(\bxi_j\right)\right) \mhE{\phi_j}{\phi_i} 
\nonumber \\
& \stackrel{\eqref{EQN::DOFI}}{=} \sum_{E \in \Th} \sum_{j = 1}^{d_{k, E}^S} 
f\left(\dofj(\uh)\right) \mhE{\phi_j}{\phi_i},
\label{EQN::INTERPOLATORY-TERM}
\end{align}
\end{linenomath}
which clearly shows that, {in the ideal case}, the {computation} of the 
nonlinear term {requires only} the matrix representation of {the bilinear 
form} $\mh{\cdot}{\cdot}$ and the evaluation of $f(\cdot)$ at the DoFs of 
{the discrete solution} $\uh$.
\begin{remark}
\label{REMMARK::STABILIZATION}
Applying stabilization in the last term of \eqref{EQN::VEM-SEMIDISCRETE} in 
the semi-discrete formulation is 
not 
necessary to obtain optimal convergence, but it can be  computationally 
convenient, as shown in Section 
\ref{SECT::FULLY-DISCRETE}.
\end{remark}
\begin{remark}
The initial condition approximation $\uh^0 = \Ih{k} u_0$ in 
\eqref{EQN::VEM-INITIAL-CONDITION} is suitable for imposing random 
initial 
data, which is commonly of interest in this kind of problems.
\end{remark}


\section{Error analysis \label{SECT::ERROR-ANALYSIS}}
This section is devoted to get an optimal error estimate in the $L_2$-norm 
for the solution of the semi-discrete formulation
\eqref{EQN::VEM-SEMIDISCRETE-FORMULATION}. The main ideas are taken from the 
error analysis carried out in \cite{Vacca_daVeiga_2015} for the {enhanced}
VEM applied to linear parabolic problems and its recent extensions to 
semilinear parabolic problems 
\cite{Adak_Natarajan_Natarajan_2019,Adak_Natarajan_Kumar_2019,
Adak_Natarajan_2020_b}. Nonetheless, in Theorem \ref{THM::ERROR-ESTIMATE} we 
address the following differences:
\begin{itemize}
	\item[$\bullet$] The approximated solution is sought in the S-VEM space  
$\VkS{k}{\Th}$.
	\item[$\bullet$] The nonlinear term is approximated by its interpolant 
$\Ih{k}{f(\uh)} \in \VkS{k}{\Th}$. The term $\mh{\Ih{k}{f(\uh)}}{\vh}$ in the 
semi-discrete formulation
\eqref{EQN::VEM-SEMIDISCRETE-FORMULATION} includes also a stabilization part, 
that was not present in the formulation in \cite{Adak_Natarajan_2020} for the 
sine-Gordon 
equation. This choice is computationally convenient when using an operator  
splitting time marching scheme, while retaining the same optimal convergence.
	\item[$\bullet$] Pure homogeneous Neumann boundary conditions are 
considered.
\end{itemize}

In what follows we will make the following assumptions on the mesh:
\begin{assumption}  
\label{ASSUMP::STAR-SHAPED}
{There exists a constant $\rho > 0$, such that every  element $E \in \Th$is 
star-shaped with respect to a ball $B:= B_{\rho h_E}(\uu{z})$ centered 
at $\uu{z} \in E$ and with radius $\rho h_E$, where $h_E := \diam{E}$. In 
addition, every edge $e$ of $E$ satisfies $|e| \geq \rho h_E$.}
\end{assumption}

The above assumption guarantees that the following condition holds: \textit{for 
each $E\in \Th$, there exists a ``virtual triangulation'' 
$\widetilde{\mathcal{T}}_E$ of $E$ such that $\widetilde{\mathcal{T}}_E$ is 
uniformly shape regular and quasi-uniform. The corresponding mesh size of the 
auxiliary triangulation $\widetilde{\mathcal{T}}_E$ is proportional to $h_E = 
\diam{E}$ and each edge of $E$ is a side of a triangle in 
$\widetilde{\mathcal{T}}_E$}.

The elliptic projection 
operator $\PPh : H^1(\Omega) \rightarrow \VkS{k}{\Th}$, is 
defined for each function $u \in H^1(\Omega)$ as the only element 
$\Ph u \in \VkS{k}{\Th}$ satisfying
\begin{equation}
\label{EQN::RITZ-OPERATOR}
\left\{
\begin{tabular}{ll}
$\ah{\Ph u}{\vh} = \ac{u}{\vh} \quad \forall \vh \in \VkS{k}{\Th}, $\\[0.1em]
$\quad \displaystyle \int_\Omega \Ph u \dx = 0. $
\end{tabular}
\right.
\end{equation}
Since $\Ph{u}$ is the solution to the variational problem 
\eqref{EQN::RITZ-OPERATOR}, by the coercivity and continuity of 
$\ah{\cdot}{\cdot}$ and the continuity of the linear functional 
$\ac{u}{\cdot}$, the projection operator $\PPh$ is well-defined. Furthermore, 
we can prove the following estimate as in \cite[Lemma 3.1]{Vacca_daVeiga_2015}.
\begin{lemma}
\label{LEMMA::ESTIMATE-RITZ}
	Let $\Omega$ be a convex domain, and $u \in H^{k + 1}(\Omega)$. Under 
Assumption \ref{ASSUMP::STAR-SHAPED}, there exists a 
constant $C_\alpha > 0$, depending on $\alpha_*$ and $\alpha^*$ in 
\eqref{EQN::STABILITY-A} but independent of $h$ such that
\begin{equation}
\label{EQN::RITZ-ESTIMATE}
\Norm{\Ph{u} - u}{L_2\left(\Omega\right)} \leq C_\alpha h^{k 
+ 1}\Seminorm{u}{H^{k + 1}\left(\Omega\right)}.
\end{equation}
\end{lemma}

Using standard arguments as in 
\cite{Beirao_Brezzi_Cangiani_Manzini_Marini_Russo_2013} and the classical 
Dupont-Scott theory in \cite{Brenner_Scott_2008}, the following estimates for 
the interpolant $\Ih{k}(\cdot)$ and the projection
$\piO{k}{\cdot}$ are obtained.
\begin{lemma}
\label{LEMMA::ESTIMATE-INTERPOLANT}
Under Assumption \ref{ASSUMP::STAR-SHAPED}, if $u \in H^{k + 1}(\Omega)$, there 
exists a positive constant $C_I$, depending 
only on 
$k$ and $\rho$, such that the interpolant $\Ih{k}{u} \in \VkS{k}{\Th}$ 
satisfies
\begin{equation}
\label{EQN::ESTIMATE-INTERPOLANT}
\Norm{u - \Ih{k}{u}}{L_2(E)} + h_E\Seminorm{u - \Ih{k}{u}}{H^1(E)} \leq C_I h_E^{k 
+ 1} \Seminorm{u}{H^{k + 1}(E)} \quad \forall E \in \Th.
\end{equation}
\end{lemma}
\begin{lemma}
\label{LEMMA::ESTIMATE-PIB} 
Under Assumptions \ref{ASSUMP::STAR-SHAPED}, for each element $E \in \Th$, if 
$u \in H^{k + 1}(E)$, there exists a 
polynomial 
$u_\pi \in \Pp{k}{E}$, and a positive constant $C_\pi$, depending only on 
$k$ and $\rho$, such that
\begin{equation}
\label{EQN::ESTIMATE-PIB}
\Norm{u - u_\pi}{L_2(E)} + h_E\Seminorm{u - u_{\pi}}{H^1(E)} \leq C_\pi h_E^{k 
+ 1} \Seminorm{u}{H^{k + 1}(E)}.
\end{equation}
\end{lemma}

In Lemma \ref{LEMMA::NORM-EQUIVALENCE}, a norm equivalence for the S-VEM space 
is introduced in order to exploit the Lipschitz property of $f(\cdot)$ {in the 
analysis}.
\begin{lemma}
	\label{LEMMA::NORM-EQUIVALENCE}
	Under Assumption \ref{ASSUMP::STAR-SHAPED}, there {exist} two positive 
constants $c_1$ and $c_2$ depending on the shape regularity and 
quasi-uniformity 
parameters of the auxiliary triangulation $\widetilde{\mathcal{T}}_E$ of $E$
such that
	\begin{equation}
	\label{EQN::NORM-EQUIVALENCE}
		c_1 h_E \Norm{\uu{\chi}(v)}{l_2} \leq \Norm{v}{L_2(E)} \leq c_2 h_E 
\Norm{\uu{\chi}(v)}{l_2},
	\end{equation}
	where $\uu{\chi}: \VkS{k}{E} \rightarrow \IR^{d_{k, E}^S}$ is defined by 
$\uu{\chi}(v) := \left(\dofi{v}\right)_{i = 1}^{d_{k, E}^S}.$
\end{lemma}

Lemma \ref{LEMMA::NORM-EQUIVALENCE} is an extension of the classical results in 
\cite{Douglas_Dupont_1975} for finite element spaces and can be proven 
following 
the arguments 
used by 
Chen and Huang in \cite[Thm. 4.5 and Corollary 4.6]{Chen_Huang_2018}. 

From the above Lemma we can derive the following important bound in our error 
analysis. 
\begin{lemma}
	\label{LEMMA::NONLINEAR-BOUND}
	The following bound holds for any $u \in \EFC{0}{\overline{E}}$ and 
${\uh} \in \VkS{k}{E}$
	\begin{equation}
	\label{EQN::NONLINEAR-BOUND}
		\Norm{\Ih{k}{f(u)} - \Ih{k}{f(\uh)}}{L_2(E)} \leq \frac{c_2 L_f}{c_1} 
{\left(\Norm{\Ih{k}{u} - u}{L_2(E)}+ \Norm{e_u}{L_2(E)}\right)},
	\end{equation}
	{where $e_u := u - \uh$.}

\end{lemma}
\proof By Lemma \ref{LEMMA::NORM-EQUIVALENCE}, and the Lipschitz continuity of 
$f(\cdot)$  we have
\begin{linenomath}
\begin{align*}
	\Norm{\Ih{k}{f(u)} - \Ih{k}{f(\uh)}}{L_2(E)} 
\stackrel{\eqref{EQN::NORM-EQUIVALENCE}}{\leq} c_2 h_E 
\Norm{\uu{\chi}\left(\Ih{k}{f(u)} - \Ih{k}{f(\uh)}\right)}{l_2} 
& 
\stackrel{\eqref{EQN::DOF-DEF}}{=} c_2 h_E \left(\sum_{i = 1}^{d_{k, E}^S} 
\left| f(u(\bxi_i)) - f(\uh(\bxi_i))\right|^2\right)^{\frac{1}{2}} \\
	& \stackrel{\eqref{EQN::LIPCHITZ-CONDITION}}{\leq} c_2 h_E L_f 
\left(\sum_{i 
= 1}^{d_{k, E}^S} \left|u(\bxi_i) - \uh(\bxi_i)\right|^2\right)^{\frac{1}{2}} \\
	& \stackrel{\eqref{EQN::DOF-DEF}}{=} c_2 h_E 
L_f\Norm{\uu{\chi}\left(\Ih{k}{u} - \uh \right)}{L_2(E)}\\
	& \stackrel{\eqref{EQN::NORM-EQUIVALENCE}}{\leq} \frac{c_2 L_f}{c_1} 
\Norm{\Ih{k}{u} - \uh }{L_2(E)},
\end{align*}
\end{linenomath}
{To conclude the proof it suffices to use the triangle inequality in the 
last 
term.} \QED

The following theorem provides the optimal error estimate for the 
semi-discrete formulation \eqref{EQN::VEM-SEMIDISCRETE-FORMULATION} under 
suitable regularity conditions for the exact solution. We will use $C$ to 
denote a generic constant independent of the mesh size $h$ and the 
arguments of the functions in the proof will be omitted unless they are 
necessary.

\begin{theorem}
\label{THM::ERROR-ESTIMATE}
Under the Assumption \ref{ASSUMP::STAR-SHAPED}. 
Let $\Omega$ be a convex domain, and $u$ and $\uh$ be the solutions of the 
variational problems 
\eqref{EQN::WEAK-FORMULATION} and \eqref{EQN::VEM-SEMIDISCRETE-FORMULATION}, 
respectively. For $u$ and $f(u)$ smooth enough, there exists a positive 
constant $C$ independent of $h$, such that for all $t \in (0, T]$ the following 
bound holds
\begin{linenomath}
\begin{align}
\nonumber
\Norm{\uh(\cdot, t) - u(\cdot, t)}{L_2(\Omega)}  \leq & C h^{k + 
1} \Big(\Seminorm{u_0}{H^{k + 1}(\Omega)} + \Norm{u_t}{L_1\left(0, t, H^{k + 
1}(\Omega)\right)} + \Norm{u_t}{L_2\left(0, t, H^{k + 
1}(\Omega)\right)}  \\
&  + \Norm{u}{L_2\left(0, 
t, H^{k + 1}(\Omega)\right)} + \Norm{f(u)}{L_2\left(0, t, H^{k + 
1}(\Omega)\right)} \Big).
\label{EQN::ERROR-ESTIMATE}
\end{align}
\end{linenomath}
\end{theorem}
\proof \ We start decomposing $e_u := u - \uh$\ as\ $e_u = \xi_u - \theta_h$, 
where 
$\xi_u 
= u 
- \Ph{u}$ and $\theta_h = \uh - \Ph{u}$. From Lemma \ref{LEMMA::ESTIMATE-RITZ} 
and the identity $u(\cdot, t) = u(\cdot, 0) + \int_{0}^t u_t(\cdot, \tau) 
d\tau$ we 
have the following bound for $\xi_u$
\begin{equation}
\label{EQN::XI-ESTIMATE}
\Norm{\xi_u(\cdot, t)}{L_2(\Omega)} \leq C h^{k + 1} 
\left(\Seminorm{u_0}{H^{k+1}(\Omega)} + \Seminorm{u_t}{L_1(0, t, 
H^{k+1}(\Omega))}\right).
\end{equation}
Therefore, in order to get the desired estimate, it only remains to bound 
$\Norm{\theta_h(\cdot, t)}{L_2(\Omega)}$. We now proceed similarly as in 
\cite{Adak_Natarajan_Kumar_2019}. Since $\theta_h \in \VkS{k}{\Th}$, adding and 
subtracting appropriate terms in 
the semi-discrete formulation \eqref{EQN::VEM-SEMIDISCRETE}, for any $\vh \in 
\VkS{k}{\Th}$ we get
\begin{linenomath}
\begin{align}
\nonumber
\mh{\dfhdt{\theta_h}}{\vh} + \ah{\theta_h}{\vh} & = -\mh{\Ih{k} f(\uh)}{\vh} - 
\mh{\dfhdt{}\Ph{u}}{\vh} - \ah{\Ph{u}}{\vh} \\
\nonumber
& \stackrel{\eqref{EQN::RITZ-OPERATOR}}{=} -\mh{\Ih{k} f(\uh)}{\vh} - 
\mh{\dfhdt{}\Ph{u}}{\vh} - \ac{u}{\vh}  \\
\nonumber
& \stackrel{\eqref{EQN::WEAK-FORMULATION}}{=} \mc{f(u)}{\vh}
- 
\mh{\Ih{k}f(\uh)}{\vh} + {\mc{\dfhdt{u}}{\vh}} - 
\mh{\dfhdt{}\Ph{u}}{\vh} \\
& \stackrel{\eqref{EQN::MH-BILINEAR}}{=} \sum_{E \in \Th} 
\Bigg[ \underbrace{\mcE{f(u)}{\vh} - \mhE{\Ih{k}f(\uh)}{\vh}}_{T_1^E} + 
\underbrace{\mcE{\dfhdt{u}}{\vh} - \mhE{\dfhdt{}\Ph{u}}{\vh}}_{T_2^E}\Bigg],
\label{EQN::GALERKIN-ORTHOGONALITY}
\end{align}
\end{linenomath}
hence, we will look for local estimates of $T_1^E$  and $T_2^E$ on each 
{element} $E \in 
\Th$.

By the $k$-polynomial consistency property 
\eqref{EQN::POLYNOMIAL-CONSISTENCY-M} of $\mh{\cdot}{\cdot}$, we can decompose 
$T_1^E$ as
\begin{linenomath}
\begin{align}
\nonumber
\mcE{f(u)}{\vh} - \mhE{\Ih{k}f(\uh)}{\vh} 
& \stackrel{\eqref{EQN::POLYNOMIAL-CONSISTENCY-M}}{=} \mcE{f(u) - 
\piO{k}{f(u)}}{\vh} \\
\nonumber
& \qquad + \mcE{\piO{k}{f(u)} - \piO{k}{\Ih{k}(f(u))}}{\vh} \\
\nonumber
& \qquad + \mhE{\piO{k}{\Ih{k}{f(u)}}  - \Ih{k}{f(\uh)} 
}{\vh} \\
& \quad = R_1 + R_2 + R_3.
\label{EQN::DECOMPOSITION-T1E}
\end{align}
\end{linenomath}
By {the} Cauchy-Schwarz inequality and Lemma \ref{LEMMA::ESTIMATE-PIB}, it 
is easy to see that 
\begin{equation}
\label{EQN::R1}
|R_1| \stackrel{\eqref{EQN::ESTIMATE-PIB}}{\leq} C_\pi h^{k + 1} 
\Seminorm{f(u)}{H^{k + 1}(E)} \Norm{\vh}{L_2(E)}.
\end{equation}
On the other hand, by the {stability} of the $L_2$-orthogonal projection 
$\piO{k}{\cdot}$ and Lemma \ref{LEMMA::ESTIMATE-INTERPOLANT} we have
\begin{equation}
\label{EQN::R2}
|R_2| \leq \Norm{f(u) - \Ih{k}f(u)}{L_2(E)} \Norm{\vh}{L_2(E)} 
\stackrel{\eqref{EQN::ESTIMATE-INTERPOLANT}}{\leq} C_I h^{k + 
1}\Seminorm{f(u)}{H^{k + 1}(E)} \Norm{\vh}{L_2(E)}.
\end{equation}
To bound $R_3$, we first observe that by {the triangle} inequality, the 
{stability} of the $L_2$-orthogonal projection and Lemmas 
\ref{LEMMA::ESTIMATE-INTERPOLANT} and \ref{LEMMA::ESTIMATE-PIB} we have
\begin{linenomath}
\begin{align}
\nonumber
	\Norm{\piO{k}{\Ih{k}{f(u)}} - \Ih{k}{f(u)} }{L_2(E)} & \leq 
	 \Norm{\piO{k}{\Ih{k}f(u)} - 
\piO{k}{f(u)}}{L_2(E)} +  \Norm{\piO{k}{f(u)} - f(u)}{L_2(E)} \\
\nonumber
& \quad + \Norm{f(u) - 
\Ih{k}{f(u)}}{L_2(E)} \\
& \leq \left(2C_I + C_\pi \right)h^{k+1} \Seminorm{f(u)}{H^{k+1}(E)}.
\label{EQN::SAV-IDENTITY}
\end{align}
\end{linenomath}

Lemma \ref{LEMMA::NONLINEAR-BOUND} and the bound \eqref{EQN::SAV-IDENTITY} 
together with the triangle inequality and the continuity of 
$\mh{\cdot}{\cdot}$ provide the following estimate for $R_3$:
\begin{linenomath}
\begin{align}
\nonumber
|R_3| & \stackrel{\eqref{EQN::STABILITY-M}}{\leq}  \beta^* 
\Norm{\piO{k}{\Ih{k}f(u)} - \Ih{k}f(\uh)}{L_2(E)} 
\Norm{\vh}{L_2(E)} \\
\nonumber
& \ \ \ \leq \beta^* \left(\Norm{\piO{k}{\Ih{k}f(u)} - 
\Ih{k}{f(u)}}{L_2(E)} + \Norm{\Ih{k} f(u) - \Ih{k}f(\uh)}{L_2(E)}\right) 
\Norm{\vh}{L_2(E)} \\
\nonumber
& \stackrel{\eqref{EQN::NONLINEAR-BOUND}}{\leq} C 
\Bigg(\Norm{\piO{k}{\Ih{k}f(u)} - 
\Ih{k}{f(u)}}{L_2(E)} + {\Norm{\Ih{k} u - 
u}{L_2(E)} + \Norm{e_u}{L_2(E)} }\Bigg)  
\Norm{\vh}{L_2(E)} \\
&  \stackrel{\eqref{EQN::SAV-IDENTITY}}{\leq} C \left(h^{k + 1} 
\Seminorm{f(u)}{H^{k + 1}(E)} + h^{k + 1} 
\Seminorm{u}{H^{k + 1}(E)} + \Norm{e_u}{L_2(E)}\right)  
\Norm{\vh}{L_2(E)}.
\label{EQN::R3}
\end{align}
\end{linenomath}
In a similar way, decomposing 
\begin{linenomath}
\begin{align*}
\mcE{u_t}{\vh} & - \mhE{\Ph{u_t}}{\vh} \\ 
& = \mcE{u_t - \piO{k}{u_t}}{\vh} - 
\mhE{\Ph{u_t} - \piO{k}{u_t}}{\vh},
\end{align*}
\end{linenomath}
and applying similar steps as before, by the commutativity of 
$\dfhdt{}\left(\cdot\right)$ and $\Ph{\cdot}${,} we get the 
following bound for 
$T_2^E$
\begin{equation}
\label{EQN::T2E-BOUND}
\left|T_2^E\right| \leq Ch^{k+1} \Seminorm{u_t}{H^{k+1}(E)}.
\end{equation}

Integrating from $0$ to $t$ at both sides of 
\eqref{EQN::GALERKIN-ORTHOGONALITY} 
and taking $\vh = \theta_h$; since $\ah{\theta_h}{\theta_h} \geq 0$, by the 
estimate 
\eqref{EQN::XI-ESTIMATE} for $\xi_u$, the bounds 
\eqref{EQN::R1}--\eqref{EQN::T2E-BOUND} and Young's 
inequality{,} we get the 
following estimate
\begin{linenomath}
\begin{align*}
\Norm{\theta_h(\cdot, t)}{L_2(\Omega)}^2 
& \stackrel{\eqref{EQN::STABILITY-M}}{\leq}  \beta_*^{-1} 
\mh{\theta_h(\cdot, t)}{\theta_h(\cdot, t)} + 2\beta_*^{-1}\int_0^t 
\ah{\theta_h}{\theta_h} {d\tau} \\
& \ \ \leq   \ 
C\mh{\theta_h(\cdot, 0)}{\theta_h(\cdot, 0)} + C\int_0^t 
\Norm{e_u}{L_2(\Omega)} \Norm{\theta_h}{L_2(\Omega)} d\tau \\
& \quad + Ch^{k + 1} \int_0^t \Bigg(
\Seminorm{f(u)}{H^{k+1}(\Omega)} + 
\Seminorm{u}{H^{k+1}(\Omega)} + \Seminorm{u_t}{H^{k+1}(\Omega)}
\Bigg)\Norm{\theta_h}{L_2(\Omega)} d\tau \\
& \ \ \leq   C\mh{\theta_h(\cdot, 0)}{\theta_h(\cdot, 0)} + C
\int_0^t \left(\Norm{\xi_u}{L_2(\Omega)} \Norm{\theta_h}{L_2(\Omega)} + 
\Norm{\theta_h}{L_2(\Omega)}^2\right) d\tau \\
 & \quad + Ch^{k + 1} \int_0^t \Bigg(
\Seminorm{f(u)}{H^{k+1}(\Omega)} + 
\Seminorm{u}{H^{k+1}(\Omega)}
+ \Seminorm{u_t}{H^{k+1}(\Omega)} \Bigg)\Norm{\theta_h}{L_2(\Omega)} d\tau \\
&\  \stackrel{\eqref{EQN::XI-ESTIMATE}}{\leq}\ 
 C \mh{\theta_h(\cdot, 
0)}{\theta_h(\cdot, 0)} + \int_0^t 
C
\Norm{\theta_h(\cdot, \tau)}{L_2(\Omega)}^2 d\tau + 
Ch^{k+1}\Bigg(\Seminorm{u_0}{H^{k+1}(\Omega)}^2 
 \\
 & \quad + \Norm{u_t}{L_1\left(0, t, 
H^{k + 
1}(\Omega)\right)}^2 + \Norm{u_t}{L_2\left(0, t, H^{k + 1}(\Omega)\right)}^2
  + \Norm{u}{L_2\left(0, 
t, H^{k + 1}(\Omega)\right)}^2   + \Norm{f(u)}{L_2\left(0, t, H^{k + 
1}(\Omega)\right)}^2 \Bigg),
\end{align*}
\end{linenomath}
and by Gr\"onwall's lemma, since $\theta_h(\cdot, 0) = \Big(\Ph{u_0} - 
u_0\Big) + \Big(u_0 - \Ih{k}{u_0}\Big)$, combined with 
the bound \eqref{EQN::XI-ESTIMATE} for $\xi_u$, and the estimates 
\eqref{EQN::RITZ-ESTIMATE}--\eqref{EQN::ESTIMATE-INTERPOLANT}, we get the 
desired estimate \eqref{EQN::ERROR-ESTIMATE} in our theorem. \QED
\begin{remark}
The term $\Norm{\piO{k}{\Ih{k}{f(u)}} - \Ih{k}{f(u)} 
}{L_2(E)}$ in \eqref{EQN::SAV-IDENTITY} must be treated with care, since a 
direct application of the bound in Lemma \ref{LEMMA::ESTIMATE-PIB} leads to the 
appearance of the undesired term $\Seminorm{\Ih{k}{f(u)}}{H^{k+1}(E)}$, that 
would become an issue in the error analysis, since the stability of the 
interpolation operator $\Ih{k}{\left(\cdot\right)}$ on the seminorm 
$\Seminorm{\cdot}{H^{k+1}(E)}$ is not guaranteed. On the other hand, bound 
\eqref{EQN::SAV-IDENTITY} is not necessary when the stability part of the last 
term of \eqref{EQN::VEM-SEMIDISCRETE} in the semi-discrete formulation is not 
considered.
\end{remark}


\section{Fully-discrete scheme\label{SECT::FULLY-DISCRETE}}
It is evident that the efficiency of any ODE solver applied to 
\eqref{EQN::VEM-SEMIDISCRETE-FORMULATION} will be greatly benefited from the 
{fast} 
evaluation of the nonlinear term in our semi-discrete formulation. In this 
paper, we choose the second order symmetric Strang operator splitting (SS-OS) 
method \cite{Strang_1968} as time marching scheme to illustrate the advantages 
of the proposed technique.

Denoting by $\bM$ and $\bA$ the matrix representation of the 
bilinear forms 
$\mh{\cdot}{\cdot}$ and $\ah{\cdot}{\cdot}$, respectively; by the identity
\eqref{EQN::INTERPOLATORY-TERM}, the semi-discrete formulation 
\eqref{EQN::VEM-SEMIDISCRETE-FORMULATION} can be written as a system of 
nonlinear differential equations as
\begin{equation}
\label{EQN::MATRIX-FORMULATION}
	\bM \dfdt{\Uh} + \bA \Uh + \bM {\uu{f}_h\left(\Uh\right)} = 
\mathbf{0},
\end{equation}
where $\Uh$ is the vector of the representation coefficients of $\uh$ in the 
basis of $\VkS{k}{\Th}$; and the {components of the vector $\uu{f}_h(\Uh)$ are 
given by $\left(\uu{f}_h\left(\Uh\right)\right)_i = \dof{i}{f(\uh)}$}. 

{In 
the 
ideal case}, $\uu{f}_h(\Uh)$ is the vector {obtained from a} component-wise 
evaluation of {the nonlinear} function 
$f(\cdot)$ at {the entries of} $\Uh$.

The SS-OS time marching scheme decomposes the system of differential equations 
\eqref{EQN::MATRIX-FORMULATION} as a series of linear and nonlinear substeps, 
usually associated to diffusion and reaction terms,  
of the form
\begin{linenomath}
\begin{subequations}
\label{EQN::SS-OS-GENERIC}
\begin{align}
\mathcal{D} & \mathcal{R}\mathcal{D} \text{ decomposition: }  
\label{EQN::SS-OS-GENERIC-DRD}
& \Uh^{(1)} = \Diff{\tau/2}{\Uh^n}, \quad  \Uh^{(2)} = 
\React{\tau}{\Uh^{(1)}}, \quad \Uh^{n + 1} = \Diff{\tau/2}{\Uh^{(2)}}, \\
\mathcal{R} & \mathcal{D}\mathcal{R}  \text{ decomposition: } 
\label{EQN::SS-OS-GENERIC-RDR}
& \Uh^{(1)} = \React{\tau/2}{\Uh^n}, \quad  \Uh^{(2)} = 
\Diff{\tau}{\Uh^{(1)}}, \quad \Uh^{n + 1} = \React{\tau/2}{\Uh^{(2)}},& 
\end{align}
\end{subequations}
\end{linenomath}
where $\tau = t_{n + 1} - t_n$ and $\Uh^n$ is the vector approximation of 
$\uh(\cdot, t_n)$.

The efficiency of combining some discontinuous Galerkin methods with an 
interpolatory approximation of the nonlinear term as spatial discretization on 
classical meshes {with} the SS-OS time marching scheme was assessed by Castillo 
and 
G\'omez in 
\cite{Castillo_Gomez_2020,Castillo_Gomez_2021}.

A necessary condition to retain the second order accuracy of 
the full SS-OS step is that each substep in \eqref{EQN::SS-OS-GENERIC} must be 
solved 
with a second order ODE 
solver itself.
{Although we are free to choose the solver for each step}, implicit methods 
might be 
more appropriate. Conversely, if an explicit method were used, we would 
face a very restrictive CFL condition {associated to} the linear substeps, 
while for 
the 
nonlinear substeps the method might 
become unstable in the case of stiff nonlinearities.

From the discussion above we decide to apply the Crank-Nicolson method to 
each substep in \eqref{EQN::SS-OS-GENERIC}. For the 
$\mathcal{D}\mathcal{R}\mathcal{D}$ decomposition 
\eqref{EQN::SS-OS-GENERIC-DRD} the resulting fully-discrete method reads
\begin{linenomath}
\begin{subequations}
\label{EQN::SS-OS-CN}
\begin{align}
	\label{EQN::DIFFUSION-1}
	\left(\bM + 
\frac{\tau}{4} \bA\right) \Uh^{(1)} & = \left(\bM - 
\frac{\tau}{4} \bA\right) \Uh^n,\\
	\label{EQN::REACTION}
\Uh^{(2)} & = 
\Uh^{(1)} - \frac{\tau}{2} \left({\uu{f}_h}\left(\Uh^{(1)}\right) + 
{\uu{f}_h}\left(\Uh^{(2)}\right)\right),  \\
\label{EQN::DIFFUSION-2}
\left(\bM + \frac{\tau}{4} \bA\right) \Uh^{n + 1} & = \left(\bM - 
\frac{\tau}{4} \bA\right) \Uh^{(2)}.  
\end{align}
\end{subequations}
\end{linenomath}
The following remarks are in order:
\begin{itemize}
	\item The linear substeps \eqref{EQN::DIFFUSION-1} and 
\eqref{EQN::DIFFUSION-2} only consist in solving two linear systems with the 
same matrix. For a fixed time step $\tau$ such matrix is even the same at any 
time, which is advantageous since a preconditioner or a full Cholesky 
factorization can be computed just once at the beginning of the 
simulation.

	\item The nonlinear substep {requires the solution of} the nonlinear 
system \eqref{EQN::REACTION}, {which is completely independent} for each 
component of the vector $\Uh^{(2)}$, and as such, 
highly parallelizable. Note that we have 
cancelled matrix $\bM$ at both sides of this equation. Such cancellation is 
only possible 
because stabilization {was also applied} to the nonlinear term in 
\eqref{EQN::VEM-SEMIDISCRETE}; {otherwise, a large coupled system of nonlinear 
equations would be obtained}. If we apply the Newton's method to 
\eqref{EQN::REACTION} each 
linear iteration $s$ reads
\begin{linenomath}
\begin{subequations}
\label{EQN::NONLINEAR-SYSTEM}
\begin{align}
\label{EQN::NONLINEAR-SYSTEM-1}
	\left(\mathbf{I} + \frac{\tau}{2}  \mathbf{D}_f\left(\Uh^{(2, 
s)}\right)\right) 
\boldsymbol{\delta}^{(s)} & = \mathbf{b}_s, \\
\Uh^{(2, s+1)} &= \Uh^{(2,s)} - \boldsymbol{\delta}^{(s)},
\end{align}
\end{subequations}
\end{linenomath}
where $\bI$ is the identity matrix, $\mathbf{b}_s = \Uh^{(2, s)} - \Uh^{(1)} + 
\frac{\tau}{2} 
\left(\uu{f}_h\left(\Uh^{(2, s)}\right) + \uu{f}_h\left(\Uh^{(1)}\right)\right) 
$ and $\mathbf{D}_f(\Uh)$ is the diagonal matrix $\mathbf{D}_f\left(\Uh\right) 
= 
\text{diag}\left(f'\left(\Uh\right)\right)$. Since 
matrix $\left(\mathbf{I} + \frac{\tau}{2} \mathbf{D}_f\left(\Uh^{(2, 
s)}\right)\right)$ is 
also diagonal, the solution of \eqref{EQN::NONLINEAR-SYSTEM-1} reduces to a 
trivial 
entry-by-entry 
division.
\end{itemize}
We end this section {with} the following well-posedness result of the 
fully-discrete scheme.
\begin{proposition} 
\label{PROP::WELL-POSEDNESS-SSOS}
The fully-discrete 
schemes $\mathcal{D}\mathcal{R}\mathcal{D}$ and 
$\mathcal{R}\mathcal{D}\mathcal{R}$ are 
well-posed for any $0 < \tau < 2/L_f$.
\end{proposition}
\PROOF \ 
Without loss of generality we will prove the well-posedness only for the 
$\mathcal{D}\mathcal{R}\mathcal{D}$ scheme. 

Since matrix 
$\left(\bM + \frac{\tau}{4} \bA\right)$ is symmetric and positive 
definite, the 
existence of the solution of each linear substep in \eqref{EQN::DIFFUSION-1} 
and \eqref{EQN::DIFFUSION-2} is guaranteed. 

On the other 
hand, each independent one dimensional problem in the nonlinear  substeps 
\eqref{EQN::REACTION} is 
equivalent to find a fixed point of the function $g(x) = a - 
\frac{\tau}{2} f(x)$ for some constant $a$, which can be {easily 
shown} to be 
a contraction as long as $0 < \tau < 2/L_f$; therefore, under such condition, 
the existence of the solution of the nonlinear substeps is also guaranteed. 

Existence and uniqueness of the full step in \eqref{EQN::SS-OS-CN} would 
then proceed from those of each susbtep.  \QED
\section{Extension to arbitrary $k$ \label{SECT::HIGH-ORDER}}
We now present an extension of the interpolatory S-VEM to the general case, 
when some internal-moment DoFs are needed. The main drawback in 
such case is that for $k \geq \eta_E $, condition \eqref{EQN::BNDRY-PROJECTION} 
is not enough to define a projection due to the existence of 
$\Pp{k}{E}$-bubbles. Hence, some additional internal-moment DoFs and a 
computable projection operator are needed. 

Since for non-convex polygons the choice of the additional DoFs is 
more involved, see \cite[Sect. 3]{DaVeiga_Brezzi_Marini_Russo_2016}, we will 
focus on the case of convex polygons.  For convex polygons, if the 
internal-moment DoFs up to order $k - \eta_E$ are added, the 
projection $\pi_{k, E}^S: \Vke{k}{E} \rightarrow \Pp{k}{E}$ defined in 
\cite{Russo_2016} as 
\begin{equation*}
	\label{EQN::BNDRY-PROJECTION-2}
\DotProd{\uu{\chi}\left(\piS{k,E}{\vh}\right)}{\uu{\chi}\left(m_{\alpha}
^E\right)}_{l_2} = 
\DotProd{\uu{\chi}\left(\vh\right)}{\uu{\chi}\left(m_{\alpha}
^E\right)}_{l_2}, \quad \alpha = 1, \ldots, r_k,
\end{equation*}
is well-defined and computable from the DoFs by definition.

For any convex polygon $E \in \Th$, if $k \geq \eta_E$, the local 
Serendipity VEM space is then defined as
\begin{equation*}
	\label{EQN::MOMENTS-SERENDIPITY-VEM}
	\VkS{k}{E} := \left\{v \in \EFC{0}{\overline{E}} :  \int_E \left(v - 
\piS{k, 
E}{v}\right) m_{\alpha}^E \dx = 0, \ r_{k - \eta_E} < \alpha \leq r_k \right\}.
\end{equation*}
and $d_{k, E}^S = kN_E + \dim(\Pp{k - 
\eta_E}{E})$. 

Unfortunately, the presence of these internal-moment DoFs prevents 
the direct extension of the variational formulation 
\eqref{EQN::VEM-SEMIDISCRETE-FORMULATION} to the case when $k$ does not satisfy 
the condition of the ideal case. This is due to the fact that the 
entries of the vector $\uu{f}_h(\Uh)$ in \eqref{EQN::MATRIX-FORMULATION} 
corresponding 
to such DoFs consist of integrals of the form 
\begin{equation*}
	\label{EQN::INTERNAL-MOMENTS-FH}
	\frac{1}{|E|}\int_E f(\uh) m_\alpha \dx, \ \alpha = 1, \ldots, r_{k - 
\eta_E},
\end{equation*}
that are not computable via the DoFs of $\uh$. To overcome this problem, we 
replace the interpolant $\Ih{k}{f(\uh)}$ 
in the semi-discrete variational formulation 
\eqref{EQN::VEM-SEMIDISCRETE-FORMULATION} by a computable quasi-interpolatory  approximation in 
the space $\VkS{k}{\Th}$ that will be denoted by $\widetilde{f}_h(\uh)$. 

For clarity, we assume that the DoFs associated to $\VkS{k}{E}$ are arranged so 
that the first $(kN_E)$ of them correspond to the boundary DoFs.
Since every function in the space $\VkS{k}{\Th}$ is uniquely determined by its 
DoFs, we set the DoFs of $\widetilde{f}_h(\uh)$ on each element $E \in \Th$ as
\begin{subequations}
\label{EQN::QUASI-INTERPOLANT}
\begin{align}
	\label{EQN::QUASI-INTERPOLANT-1}
	\dof{i}{\widetilde{f}_h(\uh)} & := \dof{i}{f(\uh)} = 
f\left(\uh\left(\bxi_i\right)\right), \quad \text{ for } i = 1, \ldots, k N_E, 
\\
	\label{EQN::QUASI-INTERPOLANT-2}
	\dof{i}{\widetilde{f}_h(\uh)} & := \dof{i}{f\left(\piO{k}{\uh}\right)} = 
\frac{1}{|E|} \int_E f\left(\piO{k}{\uh}\right) m_{\alpha(i)}^E 
\dx, \  \text{ for } i = kN_E + 1, \ldots, d_{k,E}^S,
\end{align}
\end{subequations}
with $\alpha(i) := i - kN_E$. Unlike the interpolant $\Ih{k}{f(\uh)}$, the 
new approximation $\widetilde{f}_h(\uh) \in \VkS{k}{E}$ is computable via the 
DoFs of $\uh$ as desired. 

\subsection{Extension of the error estimate}
Most steps in the proof of the error estimate in Theorem 
\ref{THM::ERROR-ESTIMATE} are still valid for this extension of the method. The 
main difference lies on the decomposition of the left-hand side of  
\eqref{EQN::DECOMPOSITION-T1E} after substituting $\Ih{k}{f(\uh)}$ by 
$\widetilde{f}_h(\uh)$, where an additional term $R_4 := \mhE{\Ih{k}{f(\uh)} - 
\widetilde{f}_h(\uh)}{\vh}$ arises. Such term can be bounded using the 
continuity of the bilinear form $\mhE{\cdot}{\cdot}$ and the following Lemma.
\begin{lemma}
\label{LEMMA::NEW-APPROXIMANT}
	Let $\left\{m_\alpha^E\right\}_{\alpha = 1}^{r_{k-\eta_E}}$ be a basis of 
$\Pp{k - \eta_E}{E}$ that is uniformly bounded in the $L_\infty$-norm 
as $\Norm{m_\alpha^E}{L_\infty(E)} \leq 1, \ \alpha = 1, \ldots, r_{k - 
\eta_E}$. For  any $\uh \in \VkS{k}{E}$, the following bound holds
	\begin{equation*}
		\label{EQN::NEW-APPROXIMANT-BOUND}
		\Norm{\Ih{k}{f(\uh)} - \widetilde{f}_h(\uh)}{L_2(E)} \leq \frac{c_2 L_f 
r_{k-\eta_E}}{c_1} \left(2\Norm{e_u}{L_2(E)} + \Norm{u - \piO{k}{u}}{L_2(E)} 
\right).
	\end{equation*}
\end{lemma}
\proof Using Lemma \ref{LEMMA::NORM-EQUIVALENCE}, the definition of 
$\widetilde{f}_h(\uh)$, and the Cauchy-Schwarz inequality we have
\begin{align*}
	\|\Ih{k}{f(\uh)} - \widetilde{f}_h(\uh) \|_{L_2(E)} 
	\stackrel{\eqref{EQN::NORM-EQUIVALENCE}}{\leq} & c_2 h_E 
\Norm{\uu{\chi}\left(\Ih{k}f(\uh) - \widetilde{f}_h(\uh)\right)}{l_2} \\
	\stackrel{\eqref{EQN::QUASI-INTERPOLANT}}{=} & \frac{c_2 h_E}{|E|} 
\left[\sum_{\alpha =1}^{ r_{k - \eta_E}} \left( \int_{E} 
\left(f\left(\uh\right) 
- f\left(\piO{k}{\uh}\right) \right) m_\alpha^E 
\dx\right)^2\right]^{\frac{1}{2}} \\ 
	\stackrel{\eqref{EQN::LIPCHITZ-CONDITION}}{\leq} & c_2 h_E L_f \Norm{\uh - 
\piO{k}{\uh}}{L_2(E)} \left[\sum_{\alpha = 1}^{r_{k - \eta_E}} 
\frac{\Norm{m_\alpha^E}{L_2(E)}^2}{|E|^2}\right]^{\frac{1}{2}} \\
	\leq \  &  \frac{c_2 L_f r_{k - \eta_E}}{c_1} \Norm{\uh - 
\piO{k}{\uh}}{L_2(E)}.
\end{align*}
The assertion follows by the triangle inequality and the stability of the 
$\piO{k}{\cdot}$ projection. \QED

An example of a polynomial basis satisfying the uniformly boundedness condition 
in the statement of the previous lemma is the scaled monomial basis defined in 
\cite{Beirao_Brezzi_Cangiani_Manzini_Marini_Russo_2013}.

\subsection{Implementation of the fully-discrete scheme}

The matrix representation of the semi-discrete variational formulation becomes 
\begin{equation*}
\label{EQN::NEW-MATRIX-FORMULATION}
	\bM \dfdt{\Uh} + \bA \Uh + \bM \widetilde{\uu{f}}_h \left(\Uh\right) = 
\mathbf{0},
\end{equation*}
where $\widetilde{\uu{f}}_h\left(\Uh\right)$ is the vector with entries given 
by 
$\left(\widetilde{\uu{f}}_h\left(\Uh\right)\right)_i := 
\dof{i}{\widetilde{f}_h(\uh)}$. Note that, as in the ideal case, the entries of 
$\widetilde{\uu{f}}_h\left(\Uh\right)$ associated to the boundary DoFs can be 
computed evaluating $f(\cdot)$ at the corresponding entries of the vector 
$\Uh$. 
As a result, the nonlinear substeps in the fully-discrete SS-OS scheme 
\eqref{EQN::SS-OS-CN} can be solved in a static condensation fashion in two consecutive steps:
\begin{enumerate}
	\item We first solve the independent one dimensional nonlinear equations 
associated to the boundary DoFs as in \eqref{EQN::NONLINEAR-SYSTEM}.
	\item Using the computed values of the boundary DoFs, we solve a set of 
independent small nonlinear systems involving just the internal-moment DoFs on 
those elements $E \in \Th$ where the condition $k < \eta_E$ is not satisfied. More specifically, for each element $E$ such that $k \geq \eta_E$, let $\UhE$ be the vector coefficient of the representation of $\uh|_E$. As the components of $\UhE$ associated to the boundary DoFs are already available from the previous step, it only remains to find the components associated to internal-moment DoFs that satisfy
\begin{equation}
\label{EQN::LOCAL-SYSTEM}
	\UhE^{(2)} = \UhE^{(1)} - \frac{\tau}{2} \left(\widetilde{\uu{f}}_h \left(\UhE^{(1)}\right) + \widetilde{\uu{f}}_h \left(\UhE^{(2)}\right)\right).
\end{equation}
We recall that by definition \eqref{EQN::QUASI-INTERPOLANT-2}, the evaluation of $\widetilde{\uu{f}}_h \left(\UhE\right)$ requires the computation of $\piO{k,E}{\uh}$ which is a local projection, i.e., it is computable using only the components of $\UhE^{(2)}$. Therefore, it is clear that the system \eqref{EQN::LOCAL-SYSTEM} is completely local and as such it can be solved separately for each element $E \in \Th$ such that $k \geq \eta_E$.
\end{enumerate}

{
\begin{remark}
\label{REMMARK::ETA-E}
The actual computation of $\eta_E$ for each $E \in \Th$ is an important and 
delicate issue in the implementation of the S-VEM. In practice, it is also 
necessary to be careful with small or almost aligned edges for stability 
reasons. We briefly recall the most used strategies in the S-VEM literature 
\cite{DaVeiga_Brezzi_Marini_Russo_2016,Beirao_Brezzi_Dassi_Marini_Russo_2018}. 
The \textit{lazy} choice consists in using always internal 
moments of degree up to $k - 3$, as by definition $\eta_E \geq 3$. A second 
strategy called the \textit{stingy} choice consists in fixing a minimum 
angle 
$\theta_0 > 0$ and then, considering as ``different'' straight lines those 
associated to consecutive edges whose internal angle is smaller than 
$\theta_0$. One last strategy, is the \textit{adaptive stingy} choice, that 
in addition to the angle treshold $\theta_0$, also impose an edge ratio 
$\rho_0$ 
and neglects those edges $e$ of $E$ satisfying $|e| < \rho_0 h_E$. Needlessly 
to say, a \textit{stingy} or \textit{adaptive stingy} choice would be more 
appropriate for the proposed method, as the additional cost of 
computing the ``exact'' value of $\eta_E$ is evidently negligible compared to 
the cost of evaluating the nonlinear term on each time step using numerical 
quadratures.
\end{remark}
}

\section{Numerical experiments\label{SECT::NUMERICAL-EXPERIMENTS}}
\begin{figure}[!ht]
\centering
\subfloat[Distorted squares mesh. ]{
\includegraphics[width=.3\textwidth,height=2.0in]{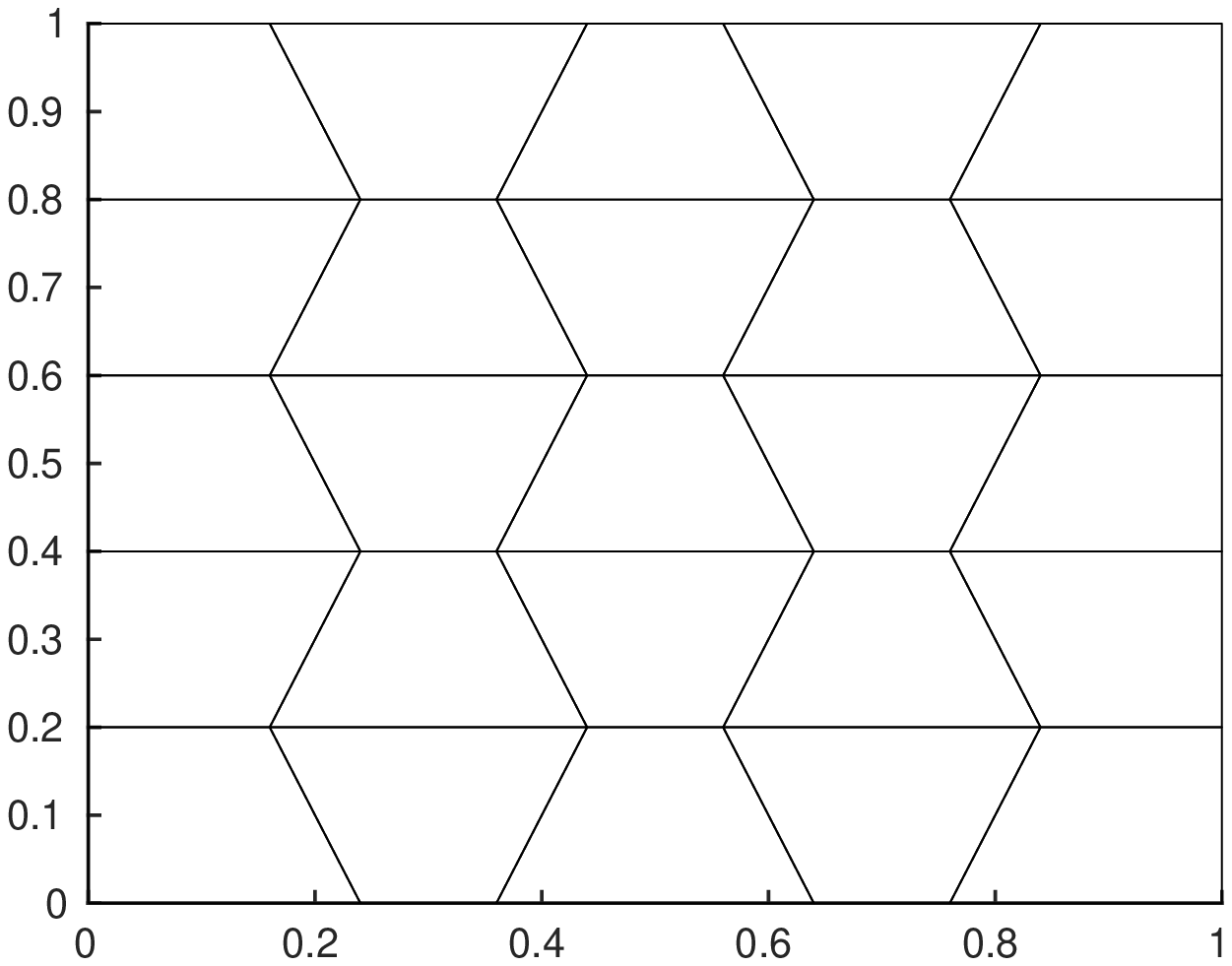}
}
\hfill
\subfloat[Voronoi mesh. ]{
\includegraphics[width=.3\textwidth,height=2.0in]{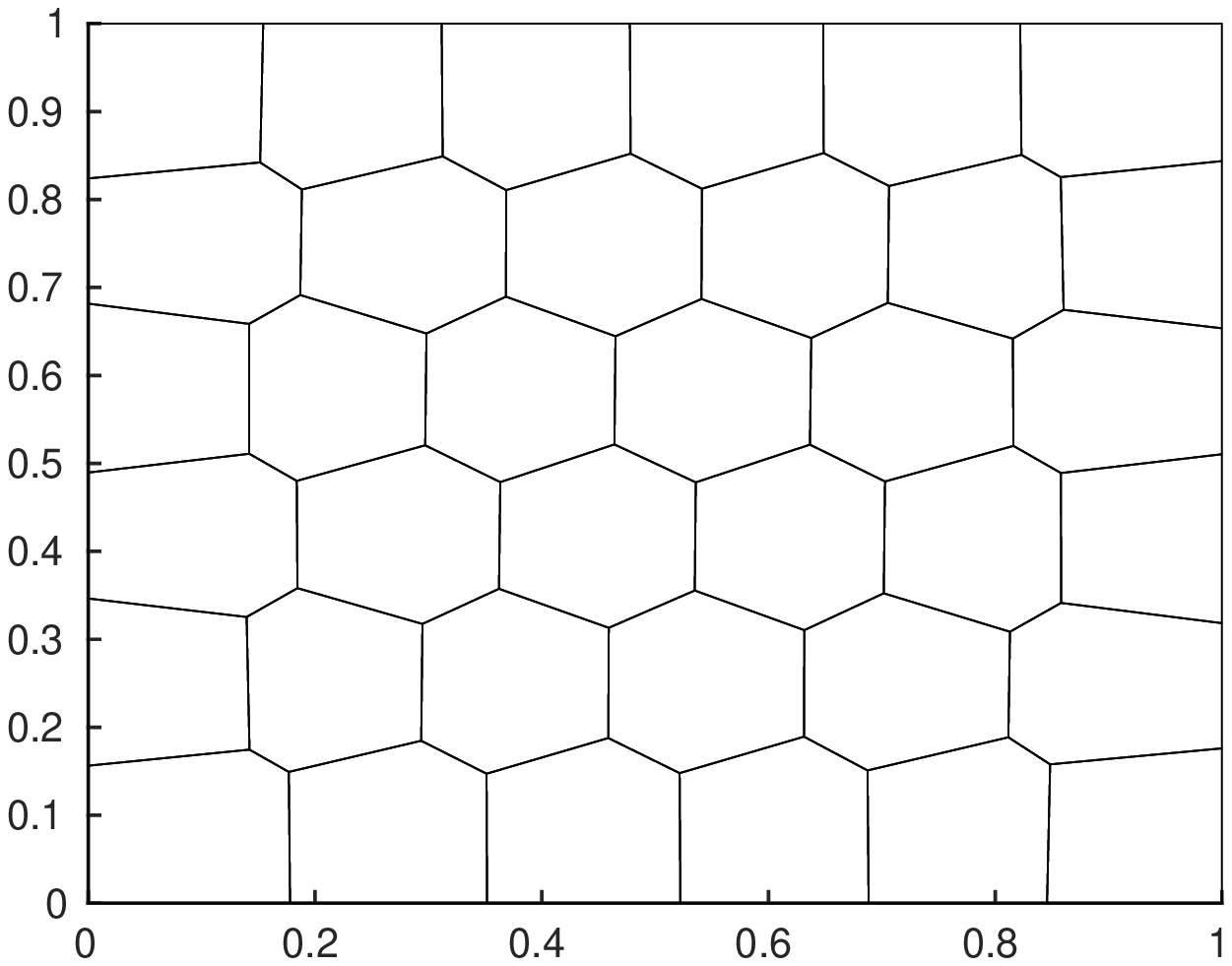} 
} 
\hfill
\subfloat[Non-convex mesh. ]{
\label{FIG::GE-MESH}
\includegraphics[width=.3\textwidth,height=2.0in]{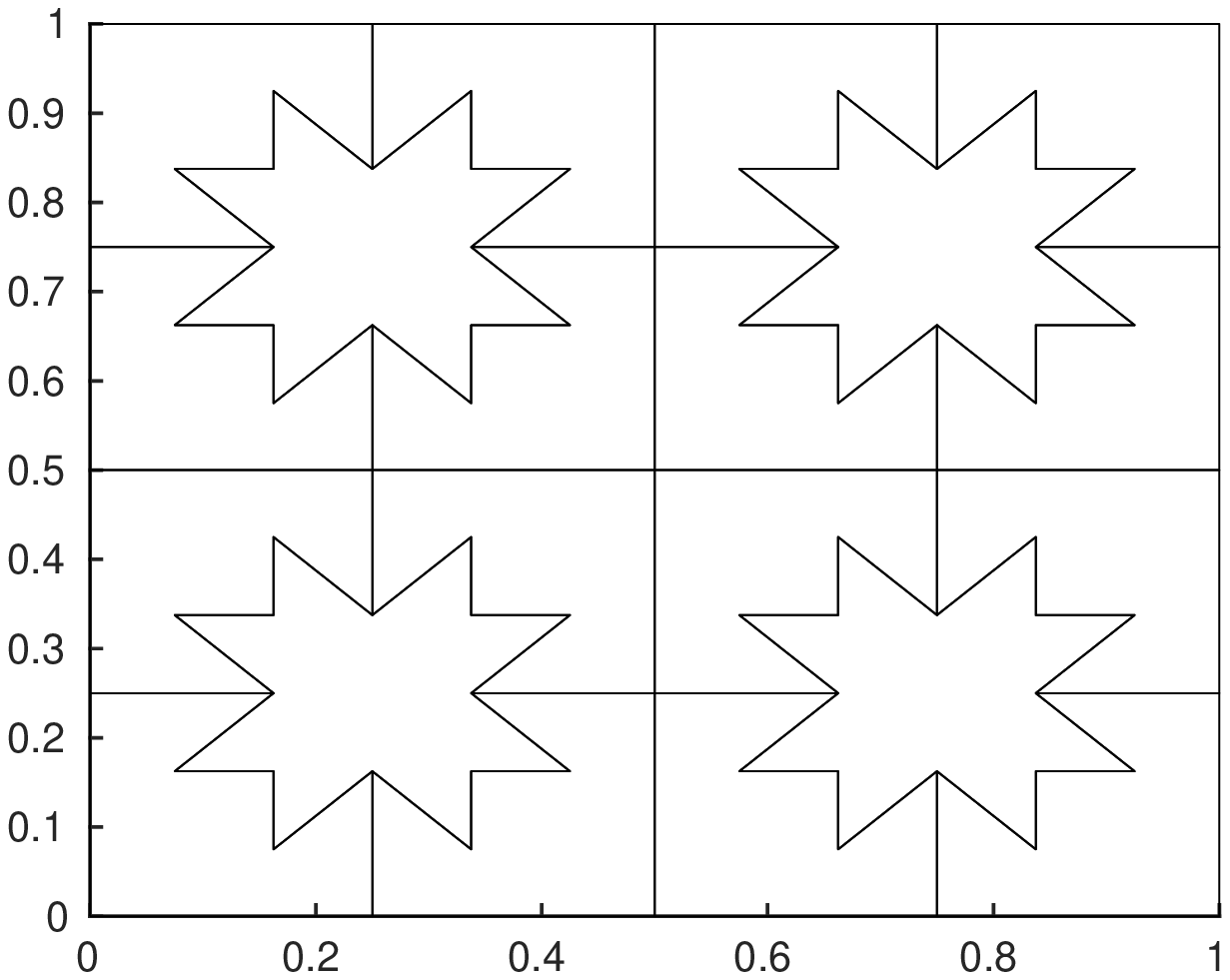}  
}
\caption{Example of meshes used in the numerical experiments.  
\label{FIG::MESHES}}
\end{figure}

In this section we present some numerical experiments to show  the accuracy and 
efficiency of the proposed scheme. An object oriented 
implementation in MATLAB was developed for high order approximations 
on general polygonal meshes. As time marching scheme we use the 
SS-OS method \eqref{EQN::SS-OS-GENERIC} presented in section 
\ref{SECT::FULLY-DISCRETE}. All the linear systems were solved {with the} 
preconditioned conjugate gradient (PCG) method. {The incomplete Cholesky 
factorization with a drop tolerance of $10^{-5}$ {was used} as preconditioner}. 
Linear and nonlinear systems were solved with a 
tolerance of $10^{-10}$ as stopping criteria; and numerical 
quadratures for each 
polygon were obtained using the Vianello approach 
\cite{Sommariva_Vianello_2007}. The sets of meshes used in all the experiments 
are exemplified in Fig. \ref{FIG::MESHES}. {Note that, for these meshes, the 
values of $\eta_E$ can be known \textit{a-priori}, {as for 
strictly convex $N$-sided polygons $\eta_E = N$ and all the non-convex polygons in 
Fig. \ref{FIG::GE-MESH} satisfy $\eta_E \geq 8$}}. 

In order to illustrate the accuracy and efficiency of the proposed method, we 
will compare 
our results with those obtained for the {enhanced VEM proposed in 
\cite{Adak_Natarajan_Kumar_2019}}. While the linear substeps of the SS-OS time 
marching scheme are similar for both versions, the nonlinear substeps for 
the method in \cite{Adak_Natarajan_Kumar_2019}
require to solve the following strongly coupled system of nonlinear 
equations
\begin{equation} 
\label{EQN::NON-INTERPOLATORY-SYSTEM}
\bM \Uh^{(2)} = \bM\Uh^{(1)} - 
\frac{\tau}{2}\left(\uu{F}_h\left(\Uh^{(1)}\right) + 
\uu{F}_h\left(\Uh^{(2)}\right)\right),
\end{equation}
where $\uu{F}_h(\cdot)$ is the nonlinear operator defined as
$${\DotProd{\uu{F}_h(\Uh)}{V_h}_{l_2} := 
\mc{\piO{k}{f(\piO{k}{\uh)}}}{\vh}} \ 
\forall \vh \in 
\VkS{k}{\Th}.$$

The nonlinear systems \eqref{EQN::NON-INTERPOLATORY-SYSTEM} will be solved 
using a semilinear iterative method, that 
avoids computing the Jacobian of the nonlinear term. Each linear iteration 
$s$ consists in solving the following linear system
\[\bM \Uh^{(2, s + 1)} = \bM\Uh^{(1)} - 
\frac{\tau}{2}\left(\uu{F}_h\left(\Uh^{(1)}\right) + 
\uu{F}_h\left(\Uh^{(2, s)}\right)\right).\]

On the other hand, the reported execution times correspond to computations 
carried out on a DELL laptop with an Intel Core i7-8750h processor, 32Gb of RAM 
 and Linux operating system.

\subsection{Accuracy test \label{SUBSECT::TEST-1}}

As first experiment we numerically asses the accuracy of the proposed method. 
We consider a manufactured problem on $Q_T = (0, 1)^2 \times (0, 1]$ with a 
nonlinear term $f(u) = 1/(1 + u^2)$, adding a source term such that the 
exact solution be $u(x, y, t) = 
e^{-t} \cos(\pi x) \cos(\pi y)$.

\begin{figure}[!ht]
\centering
\subfloat[Distorted squares meshes. ]{
\includegraphics[width=.32\textwidth,height=2.0in]{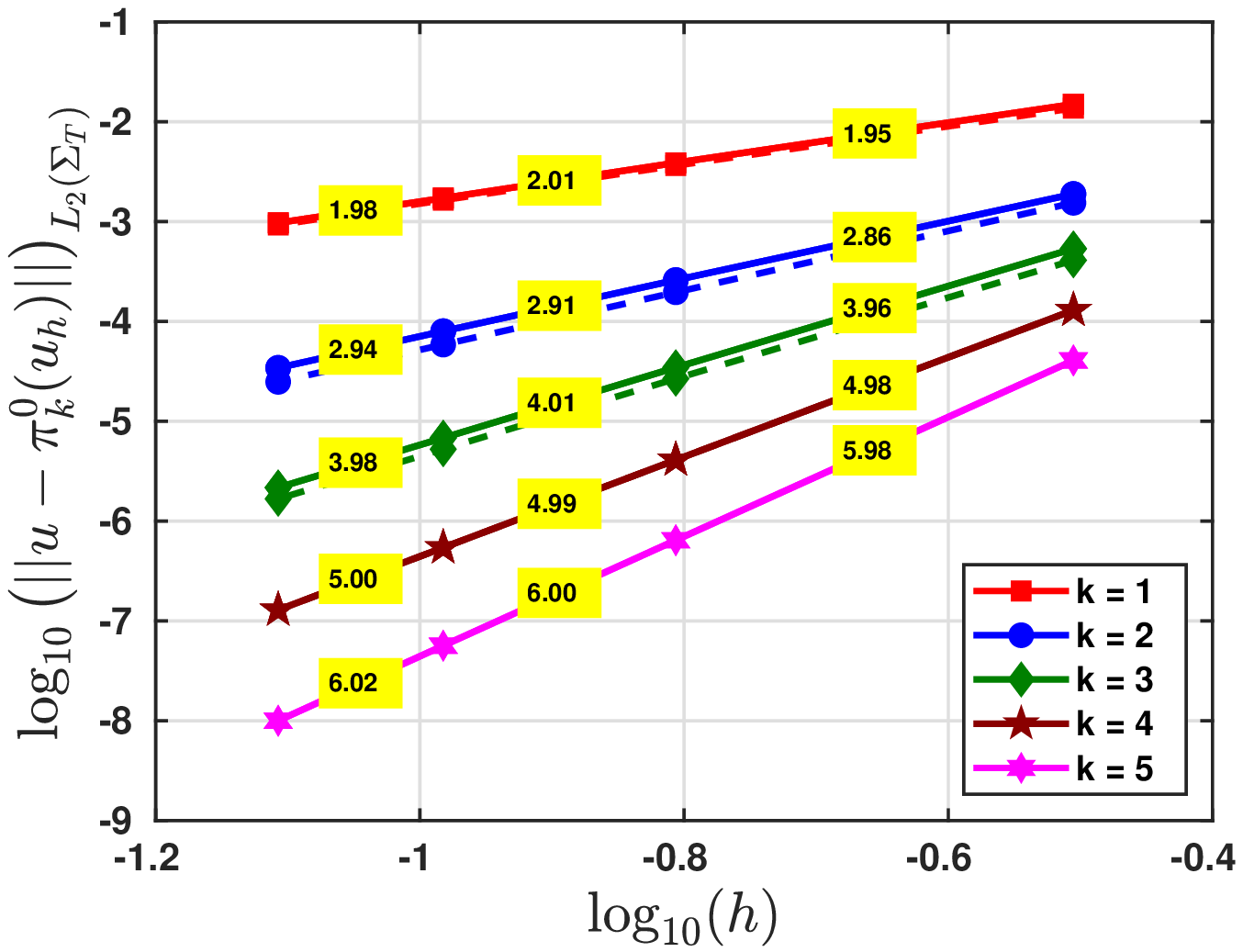}
}
\subfloat[Voronoi meshes. ]{
\includegraphics[width=.32\textwidth,height=2.0in]{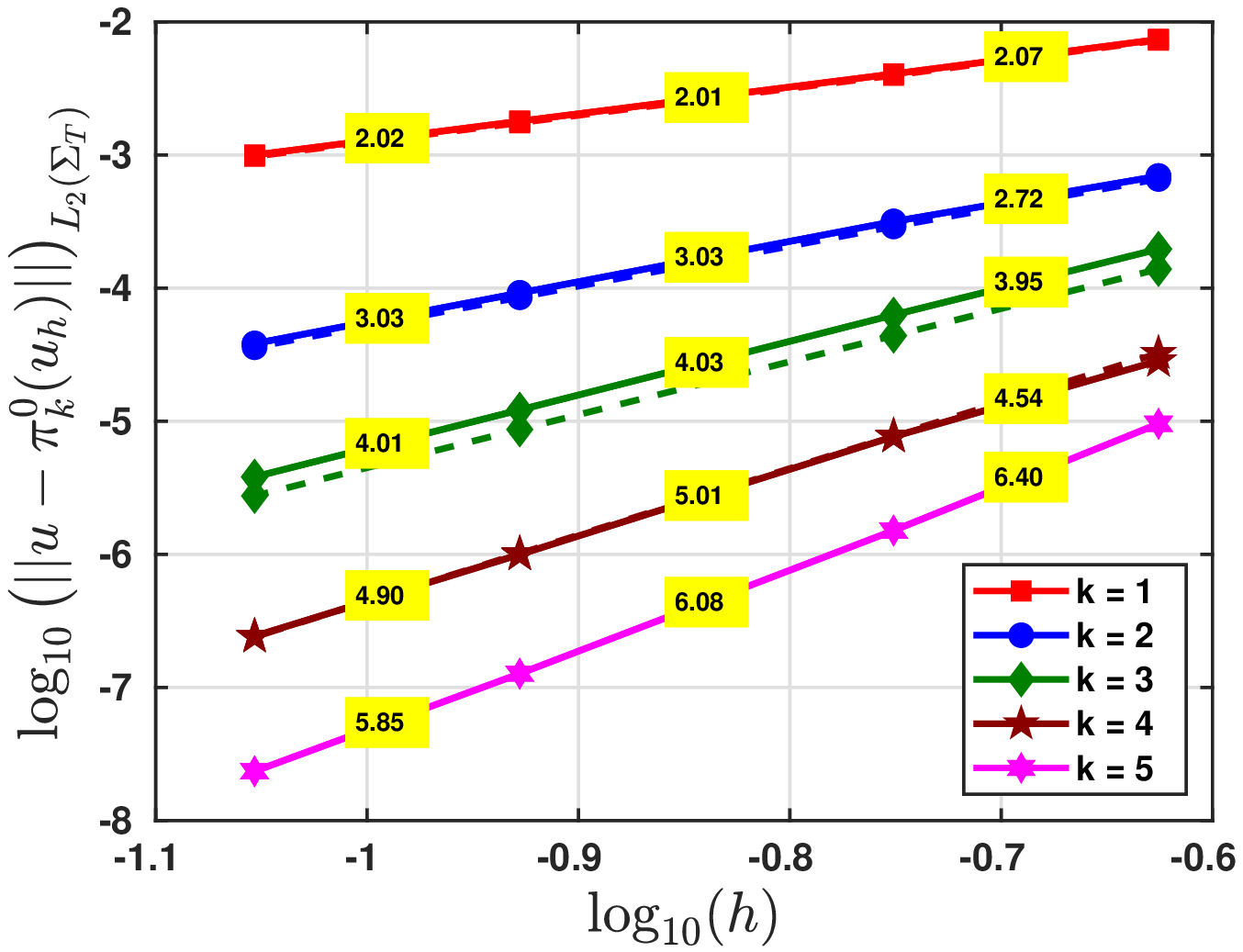}
}
\subfloat[Non-convex meshes. ]{
\includegraphics[width=.32\textwidth,height=2.0in]{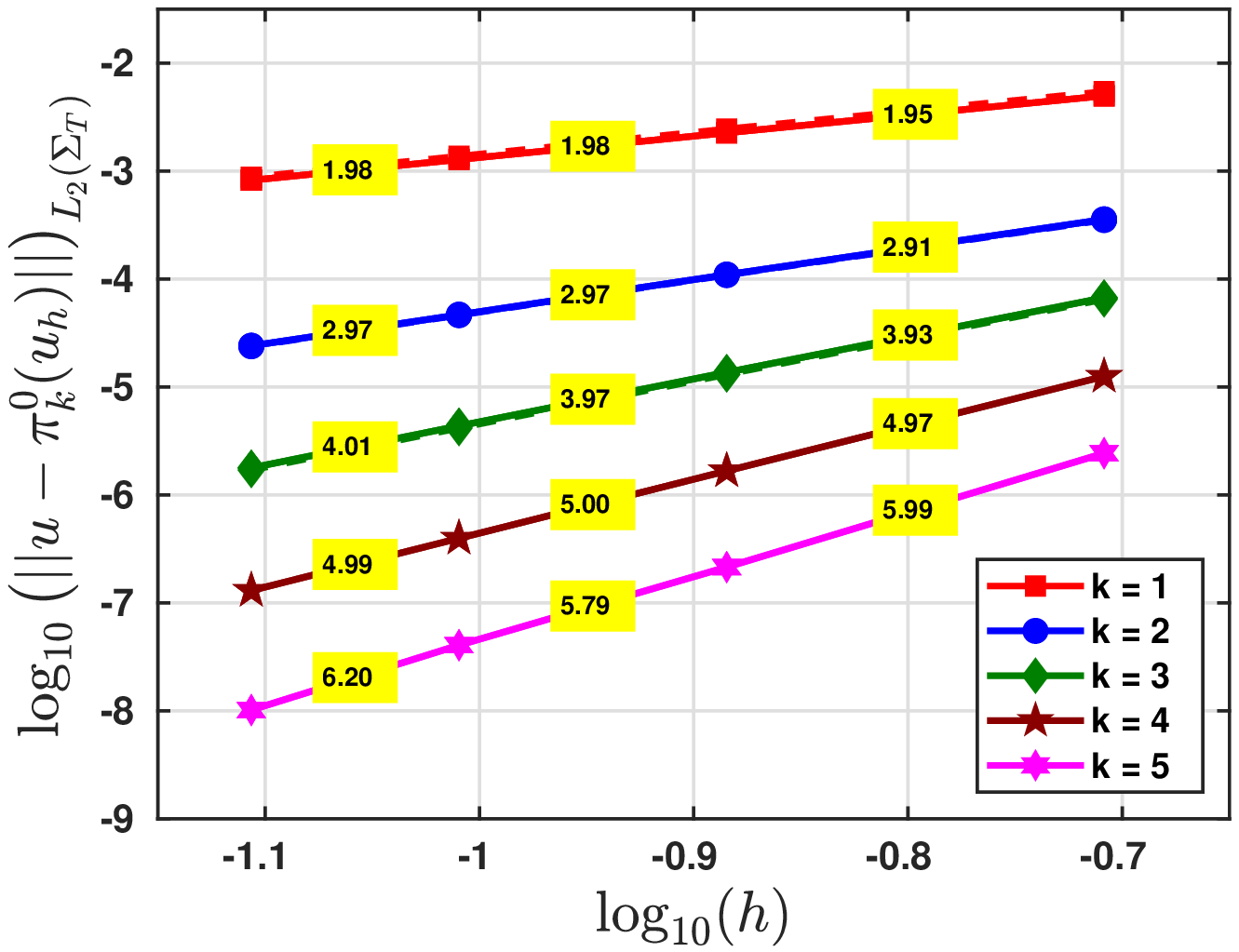}
}
\caption{Rates of convergence at $T = 1$ {for the test problem 
\ref{SUBSECT::TEST-1}} obtained for the proposed method (continuous line) and 
{the enhanced VEM in \cite{Adak_Natarajan_Kumar_2019}} (dashed line). The 
numbers in the yellow rectangles are 
the algebraic convergence rates in $h$ and {non-visible lines} were overlapped.
\label{FIG::SPATIAL}}
\end{figure}

In Fig. \ref{FIG::SPATIAL} we present the errors in the $L_2$-norm {with 
respect 
to} 
$\piO{k}{\uh}$ at the final time $T$, i.e., at $\Sigma_T := \Omega \times 
\left\{T\right\}$, for each kind of mesh. In the same plot we have included 
the errors 
obtained {by the enhanced VEM in \cite{Adak_Natarajan_Kumar_2019}} as 
reference; and no significant {difference in terms of accuracy} is observed. 
Optimal rates of 
convergence 
of order 
$\ORDER{h^{k + 1}}$ are obtained as stated in Theorem 
\ref{THM::ERROR-ESTIMATE}. The time step was 
taken as $\tau = \ORDER{h^{(k + 1)/2}}$ in order to equilibrate the errors in 
space and 
time.

\begin{figure}[!ht]
\centering
\subfloat[$\mathcal{D}\mathcal{R}\mathcal{D}$. ]{
 \includegraphics[width = .32\textwidth, height = 2.0in]{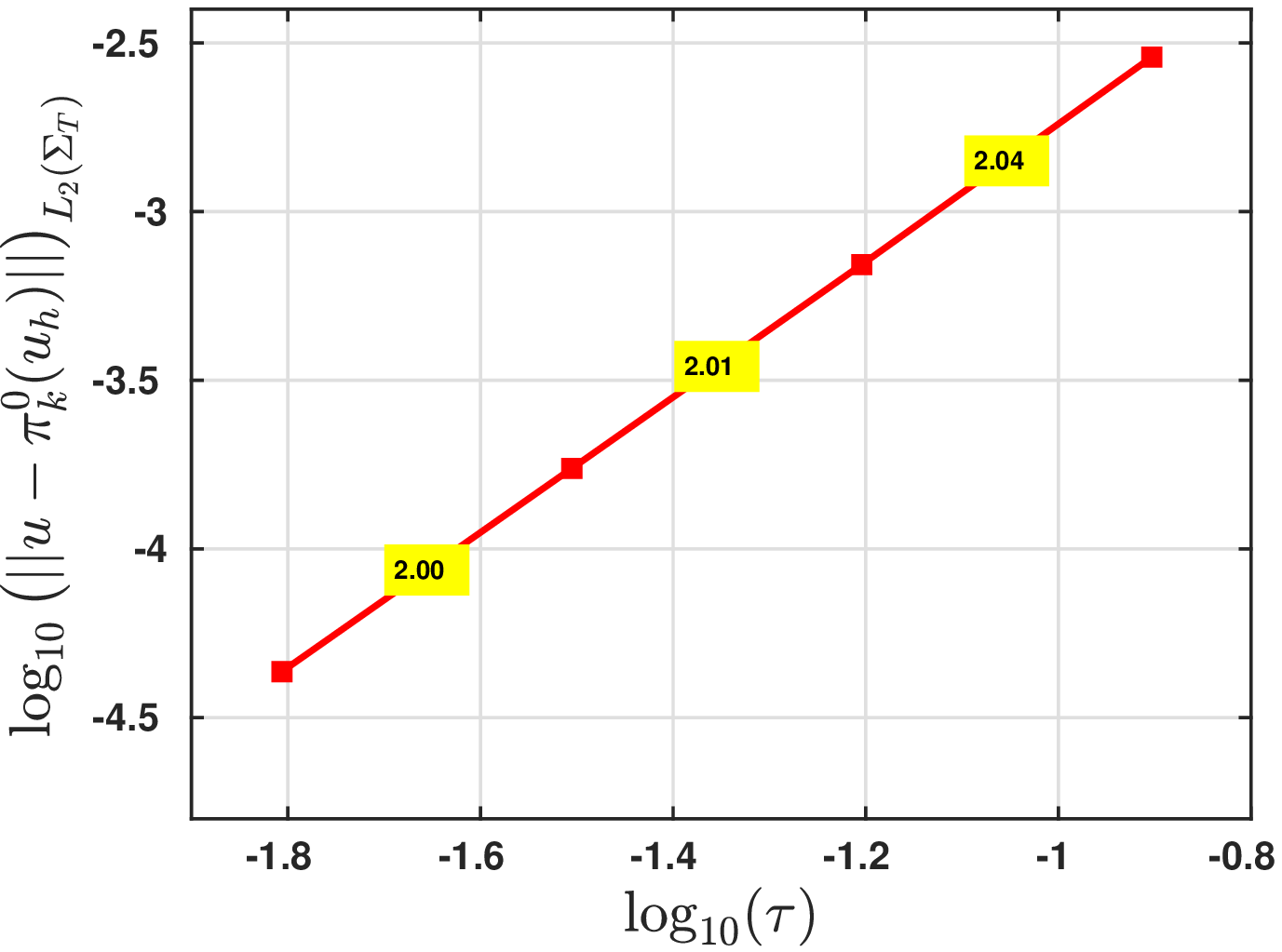}
}
\hspace{0.5cm}
\subfloat[$\mathcal{R}\mathcal{D}\mathcal{R}$. ]{
 \includegraphics[width = .32\textwidth, height = 2.0in]{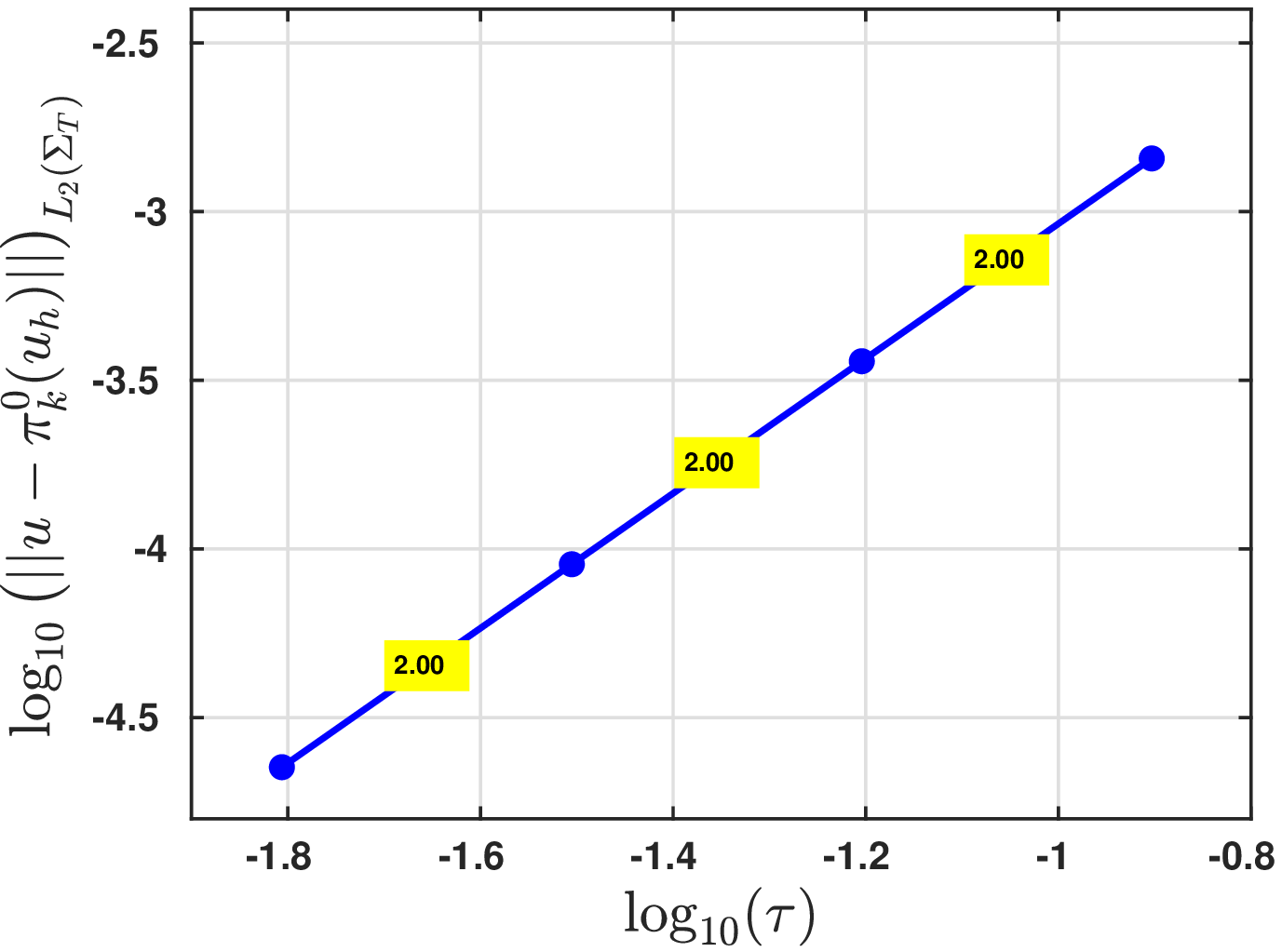}
}
\caption{Time accuracy of the proposed method at $T = 1$ {for the test 
problem \ref{SUBSECT::TEST-1}} of both versions of 
the fully-discrete scheme \eqref{EQN::SS-OS-GENERIC}. 
\label{FIG::TIME}}
\end{figure}

To evaluate the temporal accuracy of the fully-discrete scheme, we use a 
sequence of time refinements with $\tau = 1.25 \times 10^{-1},\ 
6.25\times 10^{-2},\ 3.125 \times 10^{-2},\ 1.5625\times 10^{-2}$; and in 
order to let the time error dominate, computations were carried out for the 
finest voronoi mesh and {$k = 4$}. 
The obtained 
rates of convergence for the $\mathcal{D}\mathcal{R}\mathcal{D}$ and the 
$\mathcal{R}\mathcal{D}\mathcal{R}$ splitting methods are {shown} in Fig. 
\ref{FIG::TIME} and validate the second order {in time} $\ORDER{\tau^2}$ 
accuracy of 
the {SS-OS fully-discrete scheme}. In this experiment, better accuracy is 
observed 
for the $\mathcal{R}\mathcal{D}\mathcal{R}$ splitting. Not 
shown here, similar results were obtained for the other meshes.
 
\begin{figure}[!ht]
\centering
\subfloat[Distorted squares meshes. ]{
\includegraphics[width = 0.32\textwidth, height = 2.0in]{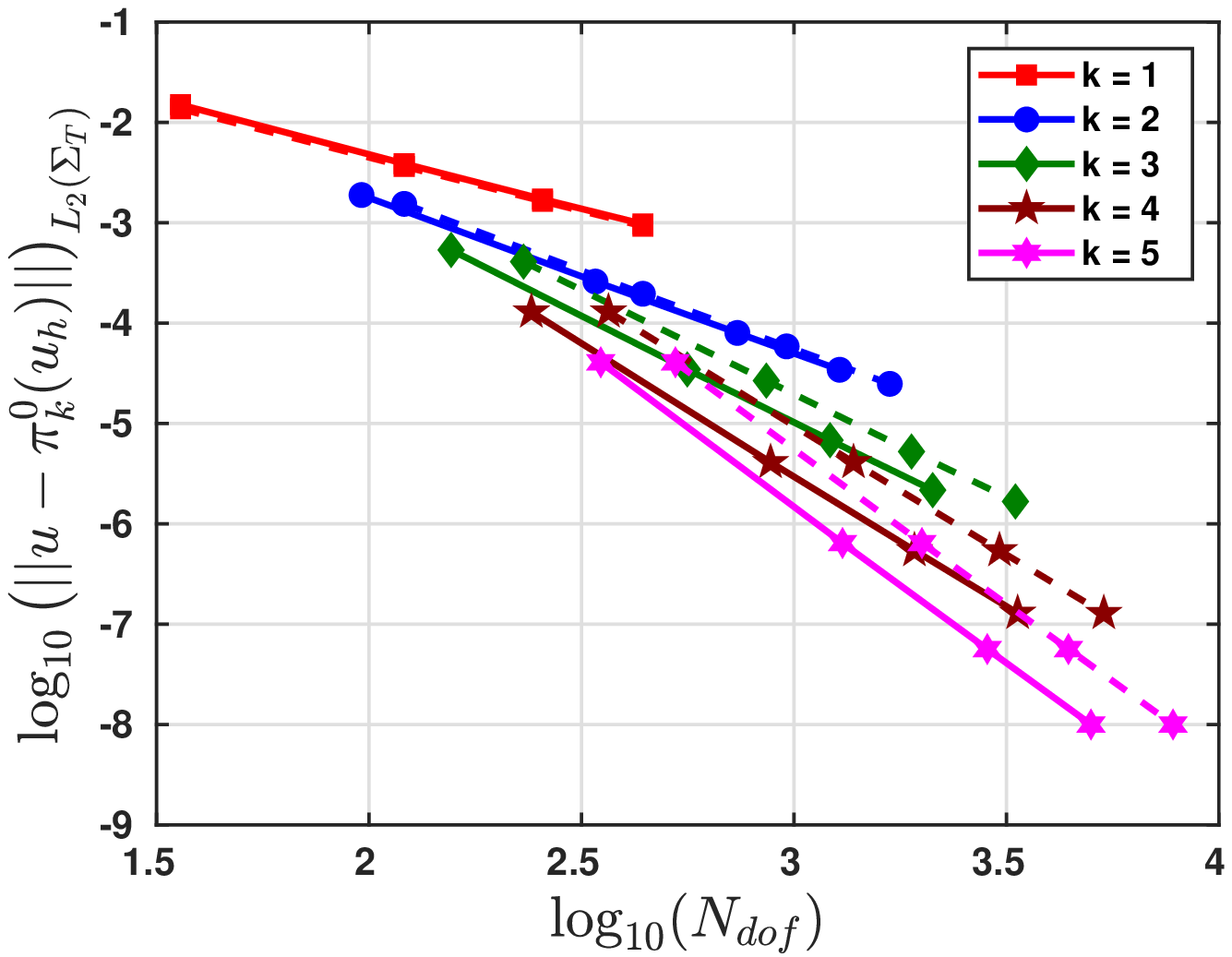}
}
\subfloat[Voronoi meshes. ]{
\includegraphics[width = 0.32\textwidth, height = 2.0in]{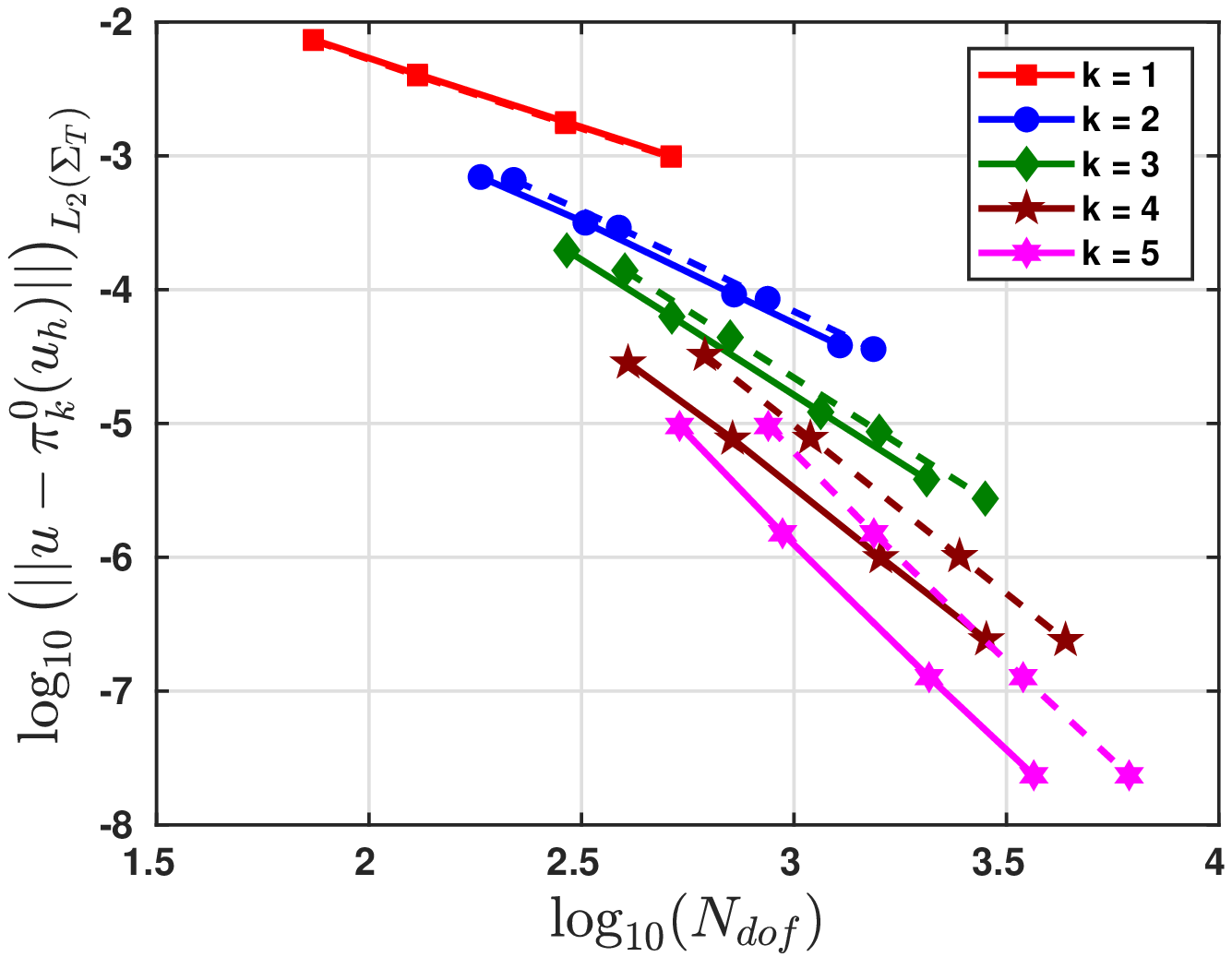}
}
\subfloat[Non-convex meshes. ]{
\includegraphics[width = 0.32\textwidth, height = 2.0in]{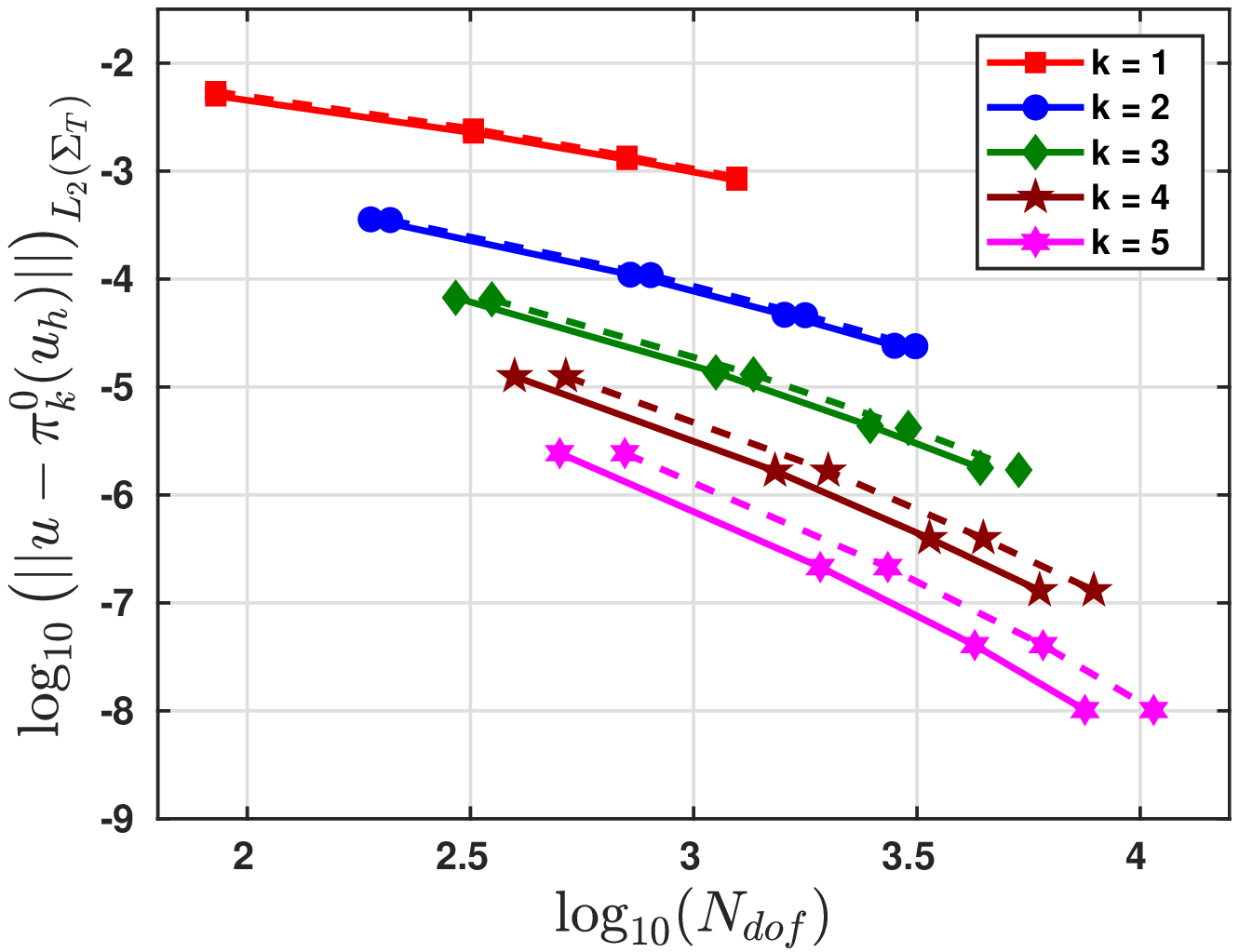}
}
\caption{{Rates of convergence at $T = 1$ {for the test problem 
\ref{SUBSECT::TEST-1} obtained by the} proposed method (continuous line) and 
{the enhanced VEM in \cite{Adak_Natarajan_Kumar_2019}} (dashed line) in 
terms of the number of DoFs.}
\label{FIG::DOFS-ERROR}}
\end{figure}

In Table \ref{TAB::N-DOF}, we compare the number of global degrees 
of freedom for the {S-VEM and the enhanced VEM in 
\cite{Adak_Natarajan_Kumar_2019}}, where naturally the reduction in the number 
of degrees of freedom depends on the mesh and a 
more noticeable reduction is obtained at increasing $k$. {This is also 
illustrated in Fig. \ref{FIG::DOFS-ERROR}, where we compare the accuracy of 
both methods with respect to the number of DoFs.}

\begin{table}
\centering
	\begin{tabular}{crrrrrr}
	Mesh type	& \multicolumn{2}{c}{Distorted squares} & 
\multicolumn{2}{c}{Voronoi} &  \multicolumn{2}{c}{Non-convex} \\
	& S-VEM & VEM \cite{Adak_Natarajan_Kumar_2019} & S-VEM & VEM 
\cite{Adak_Natarajan_Kumar_2019} & S-VEM & VEM 
\cite{Adak_Natarajan_Kumar_2019}\\
	\hline\hline
		\multicolumn{7}{c}{$k = 2$} \\
	\hline \hline
	$m_1$ & \bf 96 & 121 & \bf 183 & 219 & \bf 189 & 209 \\
	$m_2$ & \bf 341 & 441 & \bf 323 & 387 & \bf 721& 801 \\
	$m_3$ & \bf736 & 961 & \bf723 & 867 & \bf 1597& 1777\\
	$m_4$ & \bf 1281 & 1681 & \bf1283 & 1539 & \bf 2817& 3137 \\
	\hline
	\hline
		\multicolumn{7}{c}{$k = 3$} \\
	\hline
	\hline
	$m_1$ & \bf156 & 231 & \bf 292 & 400 & \bf 293 & 353 \\
	$m_2$ & \bf561 & 861 & \bf 516 & 708 & \bf 1121 & 1361 \\
	$m_3$ & \bf1216 & 1891 & \bf1156 & 1588 & \bf 2485& 3025\\
	$m_4$ & \bf 2121 & 3321 & \bf2052 & 2820 & \bf 4385 & 5345 \\
	\hline
	\hline
		\multicolumn{7}{c}{$k = 4$} \\
	\hline
	\hline
	$m_1$ & \bf \textcolor{red}{$\star$}\hspace{0.16cm} 241 & 366 & 
\textcolor{red}{$\star$}\hspace{0.16cm} \bf 407 & 617 & \bf 397& 517 \\
	$m_2$ & \bf \textcolor{red}{$\star$}\hspace{0.16cm} 881 & 1381 & 
\textcolor{red}{$\star$}\hspace{0.16cm} \bf 717 & 1093 & \bf 1521 & 2001 \\
	$m_3$ & \bf \textcolor{red}{$\star$} 1921 & 3046 & 
\textcolor{red}{$\star$} \bf 1601 & 2453 & \bf 3373 & 4453 \\
	$m_4$ & \bf \textcolor{red}{$\star$} 3361 & 5361 & 
\textcolor{red}{$\star$} \bf 2837 & 4357 & \bf 5953 & 7873 \\
	\hline
	\hline
		\multicolumn{7}{c}{$k = 5$} \\
	\hline
	\hline
	$m_1$ & \bf \textcolor{red}{$\star$}\hspace{0.16cm} 351 & 526 & 
\textcolor{red}{$\star$} \hspace{0.16cm} \bf 538 & 870  & \bf 501 & 701 \\
	$m_2$ & \bf \textcolor{red}{$\star$} 1301 & 2001 & 
\textcolor{red}{$\star$} \hspace{0.16cm} \bf 940 & 1542  & \bf 
1921 & 2721 \\
	$m_3$ & \bf \textcolor{red}{$\star$} 2851 & 4426 & 
\textcolor{red}{$\star$}\hspace{0.08cm} \bf 2080 & 
3462  & \bf 
4261 & 6061  \\
	$m_4$ & \bf \textcolor{red}{$\star$} 5001 & 7801 & 
\textcolor{red}{$\star$}\hspace{0.18cm}\bf 3668 & 6150   & \bf 
7521  & 10721 \\
	\hline 
	\hline
	\end{tabular}
	\caption{Comparison in terms of the number of degrees of freedom of the 
S-VEM and {the enhanced VEM in \cite{Adak_Natarajan_Kumar_2019} for the test 
problem \ref{SUBSECT::TEST-1}}. {The red star symbol 
$\textcolor{red}{\star}$ 
indicates that the extended version from Section \ref{SECT::HIGH-ORDER} is
needed}. \label{TAB::N-DOF} }
\end{table}

\subsection{Efficiency test \label{SUBSECT::EFFICIENCY}}

In this experiment we consider the following Allen-Cahn equation on $Q_T = (0, 
1)^2 \times (0, 22.5]$ as in \cite{Adak_Natarajan_Kumar_2019}:
\begin{linenomath}
\begin{subequations}
\label{EQN::ALLEN-CAHN}
\begin{align}
\dfhdt{u} - \epsilon \Delta u + u^3 - u & = 0, \ \hspace{3.1cm}
\mbox{ in } Q_T,\\
\nabla u \cdot \nv &= 0, \quad\quad\ \hspace{2.4cm} \mbox{ on } \partial \Omega 
\times (0, 
T),\\
u(x, y, 0) & = \cos(2\pi x^2)\cos(2\pi y^2), \hspace{0.4cm} \mbox{ in }\Omega,
\end{align}
\end{subequations}
\end{linenomath}
where the nonlinear term $f(u) = u^3 - u$ only satisfies a local Lipschitz 
condition. {In fact, the error estimate in Theorem \ref{THM::ERROR-ESTIMATE} 
is still valid if $f(\cdot)$ is only locally Lipschitz continuous 
under the additional assumption of both the exact and the approximated 
solutions 
to 
be bounded.}

{
In order to show the efficiency of the proposed method, we compare our results 
with those obtained for the interpolatory VEM in \cite{Adak_Natarajan_2020} and 
the enhanced VEM in \cite{Adak_Natarajan_Kumar_2019}. In all these experiments, 
we consider the finest meshes of each kind, $\tau = 5 \times 10^{-3}$ as 
time step and the $\mathcal{R}\mathcal{D}\mathcal{R}$ splitting.
}

{
In Table \ref{TAB::EXECUTION-ALLEN-CAHN-K1} we report the CPU execution times 
for the approximation of the Allen-Cahn equation \eqref{EQN::ALLEN-CAHN} with 
$\epsilon = 0.01$  for the proposed method and the interpolatory VEM 
presented in \cite{Adak_Natarajan_2020}. We recall that the method in 
\cite{Adak_Natarajan_2020} is limited to $k = 1$ and does not 
include the stability part of the nonlinear term, so the nonlinear systems in 
the SS-OS fully-discrete scheme \eqref{EQN::SS-OS-GENERIC} remain coupled. 
We observe that the times in the linear substeps are approximately 
equal in both cases, which is expected as both methods have the same number 
of DoFs. However, for the nonlinear substeps our method 
performs about 20 to 70  times faster depending on the mesh; and a total boost 
of approximately 10 times is obtained in all the cases.}

\begin{table}[!ht]
{
\centering
\begin{tabular}{ccccccccc}
\multicolumn{3}{c}{\bf Linear substeps} & \multicolumn{3}{c}{\bf Nonlinear
substeps} & \multicolumn{3}{c}{\bf Total} \\
S-VEM &VEM 
\cite{Adak_Natarajan_2020} & 
Ratio & S-VEM & VEM 
\cite{Adak_Natarajan_2020}& 
Ratio & S-VEM & VEM 
\cite{Adak_Natarajan_2020} & 
Ratio \\
 (sec) & (sec) & 
& (sec) & (sec)& 
& (sec) &  (sec) \\
\hline
\hline
\multicolumn{9}{c}{Distorted squares mesh} \\
2.8 & 2.6 & \bf 0.9 & 1.6 & 39.7 & \bf 24.8 & 4.4 & 42.3  & \bf 9.6 \\
\hline
\hline
\multicolumn{9}{c}{Voronoi mesh} \\
3.1 & 2.9 & \bf 0.9 & 1.5 & 48.8 & \bf 32.5 & 4.6 & 51.7 & \bf 11.2\\
\hline
\hline
\multicolumn{9}{c}{Non-convex mesh} \\
 6.6 & 6.6 & \bf 1.0 & 2.3 & 163.5 & \bf 71.1 & 8.9 & 107.1 & 
\bf 12.0\\
\hline
\hline
\end{tabular}
\caption{{CPU execution times for the Allen-Cahn equation 
\eqref{EQN::ALLEN-CAHN} {in the test problem \ref{SUBSECT::EFFICIENCY}} with 
$\epsilon = 0.01$, for the proposed S-VEM and the interpolatory VEM in 
\cite{Adak_Natarajan_2020} with $k = 1$. 
\label{TAB::EXECUTION-ALLEN-CAHN-K1}}}
}
\end{table}

{In a similar way, in Table \ref{TAB::EXECUTION-ALLEN-CAHN} we compare the 
CPU execution times for the proposed method and the enhanced VEM in 
\cite{Adak_Natarajan_Kumar_2019} with different degrees of approximation. Since 
the proposed method requires less DoFs, it performs faster for the linear 
substeps. As for the nonlinear substeps, our method performs from 40 to 
2500 times faster depending on the mesh and the degree of accuracy. A 
total boost of about 12 to 110 times is obtained. For each 
mesh, we have indicated those degrees where some internal moment DoFs are 
needed; in such cases, the extended version from Section 
\ref{SECT::HIGH-ORDER} was used and a significant improvement in the 
efficiency of the method is still observed. The substantial reduction obtained 
for the non-convex mesh is a consequence of the high number of quadrature 
points required for the VEM in \cite{Adak_Natarajan_Kumar_2019}  to compute 
the nonlinear term on each time step.}

\begin{table}[!ht]
\centering
{
\begin{tabular}{rrrrrrrrrr}
\hline
\multirow{3}{*}{ \centering $k$} & 
\multicolumn{3}{c}{ \bf \small{Linear substeps}} & 
\multicolumn{3}{c}{ \bf \small{Nonlinear
substeps}} & \multicolumn{3}{c}{\bf \small{Total}} \\
& {S-VEM} &VEM 
\cite{Adak_Natarajan_Kumar_2019} & 
Ratio & S-VEM & VEM 
\cite{Adak_Natarajan_Kumar_2019}& 
Ratio & S-VEM & VEM 
\cite{Adak_Natarajan_Kumar_2019} & 
Ratio \\
& (sec) & (sec) & 
& (sec) & (sec)& 
& (sec) &  (sec) \\
\hline
\hline
\multicolumn{10}{c}{Distorted squares mesh} \\
\hline
\hline
\hfill \raggedright  2 & 16.5 & 22.1 & \bf 1.3 & 3.6 & 374.8 & \bf 104.1 & 
20.1 & 396.9 & \bf 
19.8 \\
\hfill \raggedright 3 & 41.5 & 107.4 & \bf 2.6 & 5.3 & 1243.4 & \bf 234.6 & 
46.8 & 1350.8 & \bf 28.9 \\
\textcolor{red}{$\star$} 4 & 48.9 & 941.0 &\bf 19.2 & 208.5 & 9483.4 &\bf 
45.5 & 257.4 &10424.4 & \bf 40.5 \\
\hline
\hline
\multicolumn{10}{c}{Voronoi mesh} \\
\hline
\hline
\hfill \raggedright 2 & 9.6 & 13.4 & \bf 1.4 & 2.2 & 210.5 &\bf 95.7  & 11.8 
& 223.9 & \bf 19.0 \\
\hfill \raggedright 3 & 30.4 & 66.0 & \bf 2.2 & 3.2 & 758.3 &\bf 237.0 & 33.6 
& 824.3 & \bf 24.5\\
\textcolor{red}{$\star$} 4 & 139.3 & 570.0  &\bf 4.1 & 20.1 & 3366.7 &\bf 
167.5 & 159.4 & 3936.7 & \bf 24.7 \\
\hline
\hline
\multicolumn{10}{c}{Non-convex mesh } \\
\hline
\hline
\hfill \raggedright 2 & 54.6 & 66.6 &\bf 1.2 & 6.3 & 687.5 &\bf 109.1 & 60.9  
& 754.1 & \bf 12.4 \\
\hfill \raggedright 3 & 111.8 & 642.3 &\bf 5.8 & 7.5  & 5426.4 &\bf 723.5
& 119.3 & 6068.7 & \bf 50.9 \\
\hfill \raggedright 4 & 197.2 & 2404.4 &\bf 12.2 & 8.3 & 21388.8 & \bf 
2577.0 & 205.5 & 23793.2 & \bf 115.8\\ 
\hline
\hline
\end{tabular}
\caption{{CPU execution times for the Allen-Cahn equation 
\eqref{EQN::ALLEN-CAHN} {in the test problem \ref{SUBSECT::EFFICIENCY}} with 
$\epsilon = 0.01$, for the proposed S-VEM and {the enhanced VEM in 
\cite{Adak_Natarajan_Kumar_2019}}. The red star symbol $\textcolor{red}{\star}$ 
indicates that the extended version from Section \ref{SECT::HIGH-ORDER} was 
needed.
\label{TAB::EXECUTION-ALLEN-CAHN}}}
}
\end{table}

\begin{figure}[!ht]
\centering
\subfloat[$t = 0.1$ ]{
 \includegraphics[width = 2.6in, height = 1.8in]{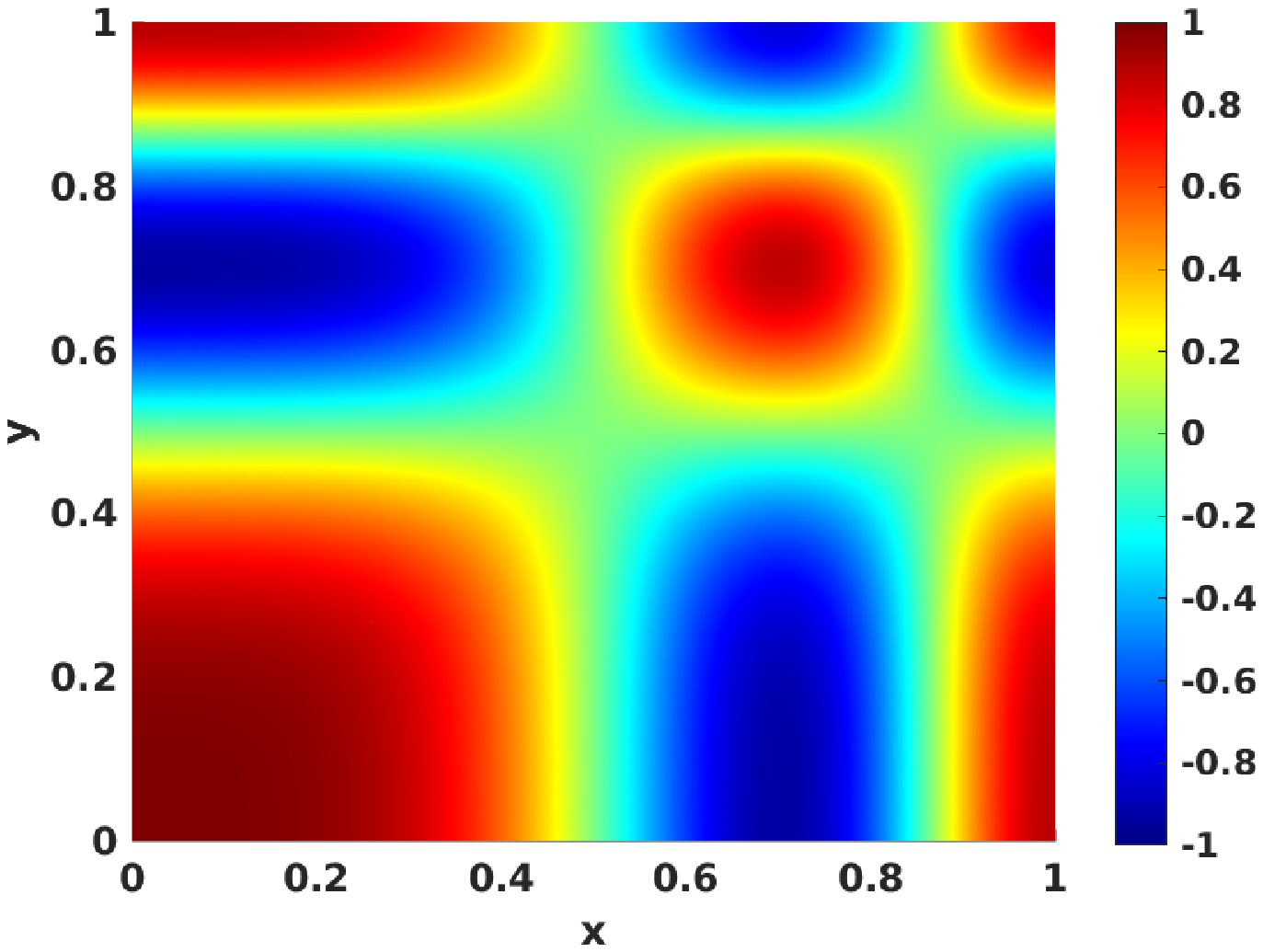}
}
\subfloat[$t = 5$ ]{
 \includegraphics[width = 2.6in, height = 1.8in]{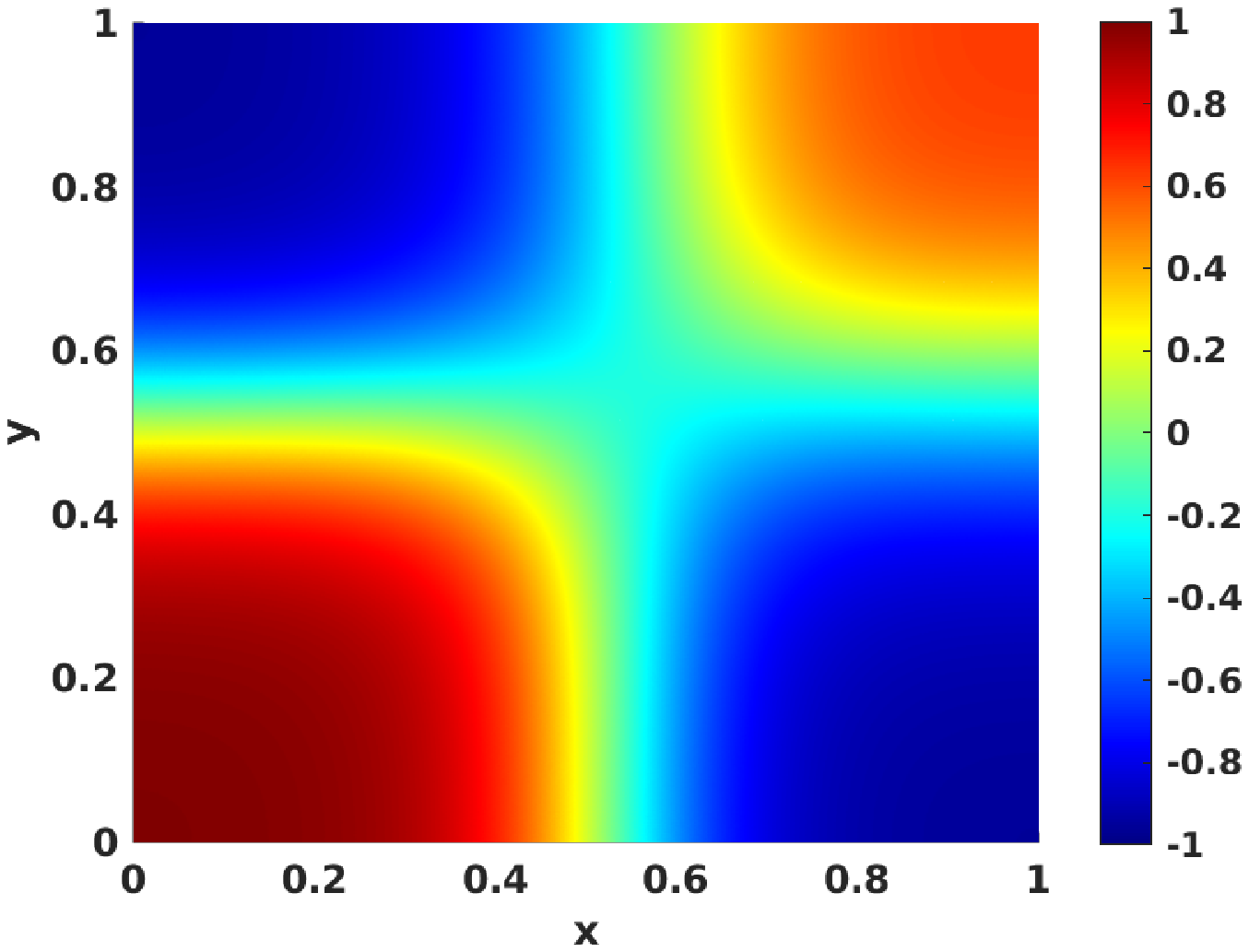}
}\\
\subfloat[$t = 10$ ]{
 \includegraphics[width = 2.6in, height = 1.8in]{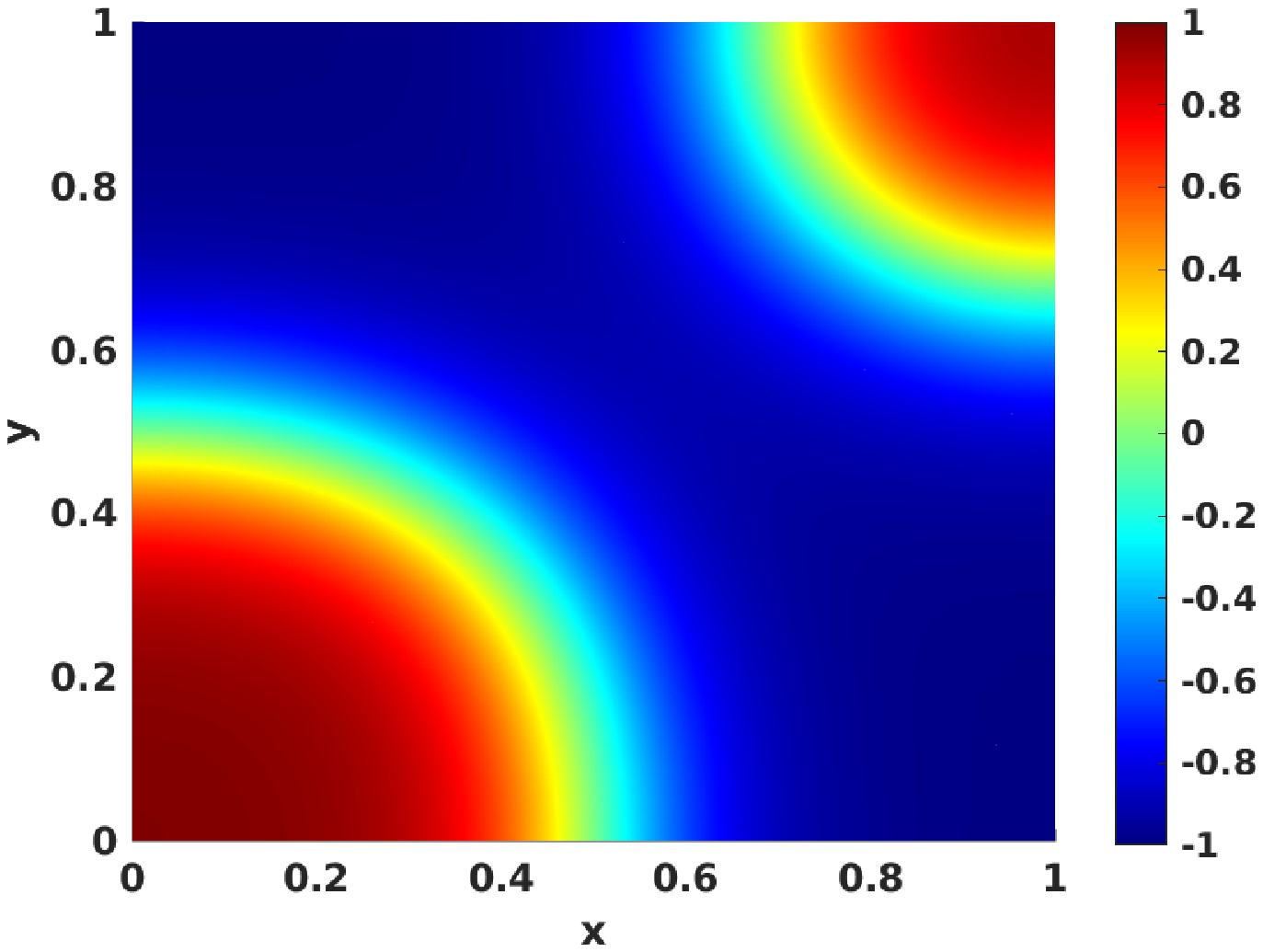}
}
\subfloat[$t = 22.5$ ]{
 \includegraphics[width = 2.6in, height = 1.8in]{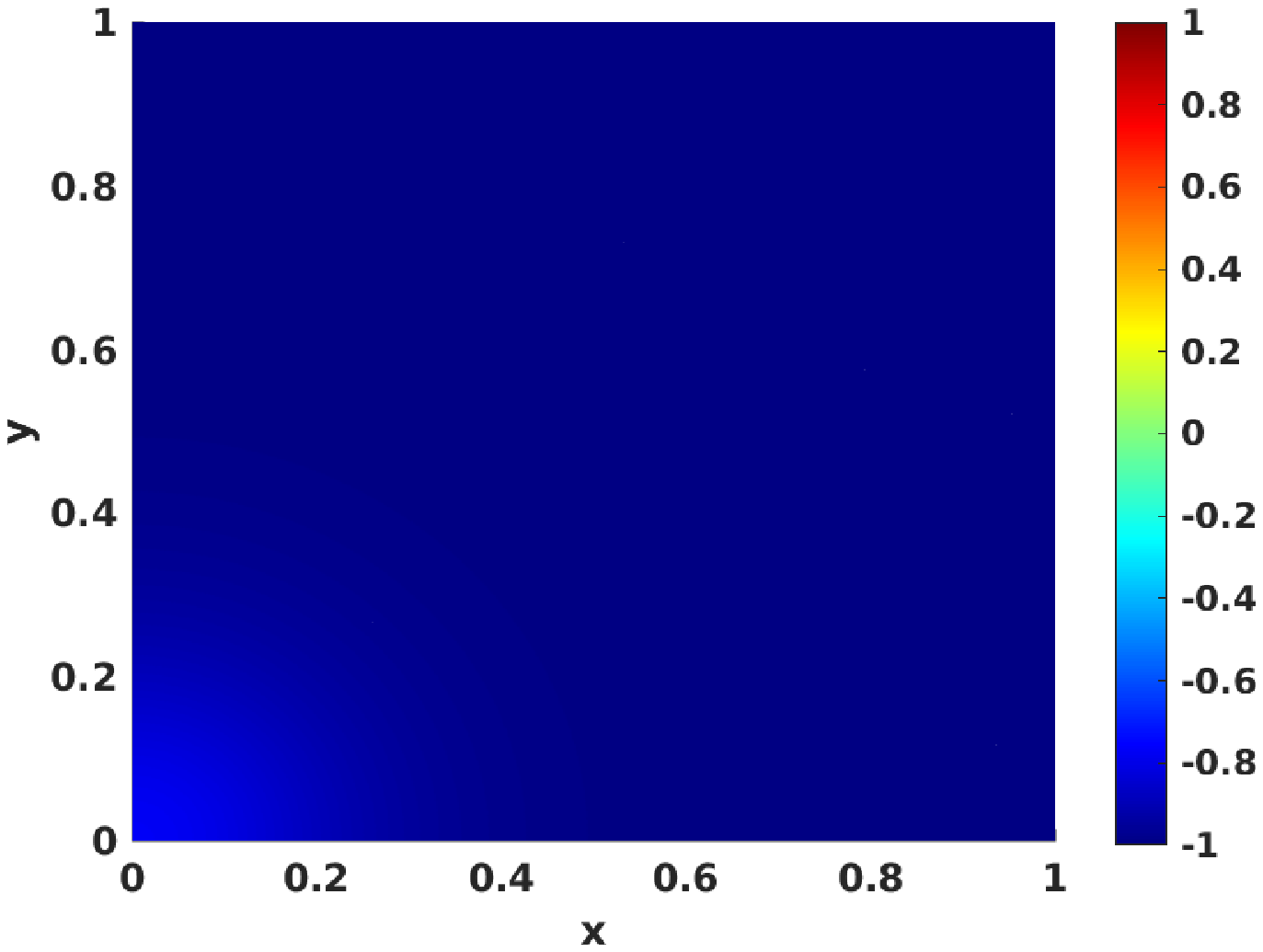}
}
\caption{Snapshots of the approximation $\piO{k}{\uh}$ for the Allen-Cahn 
equation 
\eqref{EQN::ALLEN-CAHN} with 
$\epsilon = 0.01$, {in the test problem \ref{SUBSECT::EFFICIENCY}}. 
\label{FIG::AC-PLOT-1}}
\end{figure}

In Figure \ref{FIG::AC-PLOT-1} we show the evolution of the approximated 
{solution} $\piO{k}{\uh}$ for the Allen-Cahn equation with 
$\epsilon =  \ 0.01$, which is expected to converge to its stable state $u = 
-1$. The plots
portray the same behaviour observed in \cite{Adak_Natarajan_Kumar_2019} 
for the enhanced VEM.

\section{Conclusions \label{SECT::CONCLUSIONS}}
In this work, an interpolatory Serendipity Virtual Element method for 
semilinear 
parabolic problems on polygonal meshes is proposed. A significant 
reduction in the computational cost of the method is obtained by approximating 
the nonlinear term with an element of the S-VEM space. Optimal error 
estimates of order $\mathcal{O}\left(h^{k+1}\right)$ in the $L_2$-norm are 
proven for the semi-discrete formulation.

To exploit the structure of the system of nonlinear differential equations 
{arising from} the semi-discrete formulation, we {use} a second order 
operator splitting {as time marching scheme}, which decouples the linear and 
nonlinear terms. 
{In the ideal case, with only boundary DoFs, the} nonlinear substeps 
{consist in solving} a set of completely independent one dimensional 
nonlinear equations{; while in the extension proposed to the case when some 
internal-moment DoFs are required, it is also necessary to solve an 
additional set of independent  small nonlinear systems on each element of the 
mesh that does not satisfy the condition of the ideal case.}
Our numerical experiments validate the optimal convergence of the method and 
the 
improvement in efficiency {respect to the enhanced VEM in 
\cite{Adak_Natarajan_Kumar_2019}}.

\end{document}